\documentclass[final]{siamltex}

\usepackage{float}
\usepackage{amsmath,amssymb}
\usepackage[pdftex]{graphicx,color}
\usepackage{subfigure}
\usepackage{multicol}
\usepackage{hyperref} 
\hypersetup{
  colorlinks=true,
  linktocpage=true,
  bookmarksnumbered=true
}
\renewcommand{\Re}{{\mathbb R}}         

\newtheorem{thm}{Theorem}[section]
\newtheorem{rem}[thm]{Remark}

\newtheorem{ex}[thm]{Example}

\newtheorem{alg}[thm]{Algorithm}
\newtheorem{defn}[thm]{Definition}

\DeclareMathOperator{\sgn}{sgn}

\newcommand{\ud}{\mathrm{d}}
\newcommand{\norm}[1]{\left\|#1\right\|}
\newcommand{\abs}[1]{\left|#1\right|}
\newcommand{\rset}{\mathbb{R}}

\title{An A Posteriori Error Estimate for Symplectic Euler Approximation of Optimal Control Problems\thanks{This work was supported by the Swedish Research Council and the Swedish e-Science Research Center. The fifth author is a member of the Research Center on Uncertainty Quantification in Computational Science and Engineering at KAUST.}}

\author{Jesper Karlsson\footnotemark[2]\ \footnotemark[3]
\and Stig Larsson\footnotemark[4]
\and Mattias Sandberg\footnotemark[5]
\and Anders Szepessy\footnotemark[5]
\and Ra\`ul Tempone\footnotemark[2]}

\begin{document}

\maketitle

\renewcommand{\thefootnote}{\fnsymbol{footnote}}

\footnotetext[2]{SRI UQ Center, Computer, Electrical, and Mathematical Sciences and Engineering, King Abdullah University of
Science and Technology (KAUST), Thuwal, Saudi Arabia ({\tt jesper.karlsson@kaust.edu.sa, raul.tempone@kaust.edu.sa})}

\footnotetext[3]{Dynamore Nordic AB, Theres Svenssons gata 10, S--417 55 G\"oteborg, Sweden ({\tt jesper@dynamore.se})}

\footnotetext[4]{Department of Mathematical Sciences, Chalmers University of
Technology and University of Gothenburg,
S--412 96 Gothenburg, Sweden ({\tt stig@chalmers.se})}

\footnotetext[5]{Department of Mathematics, KTH Royal Institute of Technology,
S--100 44 Stockholm, Sweden ({\tt msandb@kth.se, szepessy@kth.se})}

\renewcommand{\thefootnote}{\arabic{footnote}}

\begin{abstract}
This work focuses on numerical solutions of optimal control problems.
A time discretization error representation is derived for the approximation of the associated
value function. It concerns Symplectic Euler solutions of the
Hamiltonian system connected with the optimal control problem. The
error representation has a leading order term consisting of an error
density that is computable from Symplectic Euler solutions. Under an
assumption of the pathwise convergence of the approximate dual function as the maximum time step goes to
zero, we prove
that the remainder is of higher order than the leading error density
part in the error representation.
With the error representation, it is possible to perform adaptive time
stepping. We apply an adaptive algorithm originally developed for
ordinary differential equations. The performance is illustrated by
numerical tests.
\end{abstract}




\numberwithin{equation}{section}
\numberwithin{figure}{section}
\numberwithin{table}{section}

\setcounter{tocdepth}{1}




\section{Introduction}
In this work, we will present an asymptotic a posteriori error estimate for optimal
control problems.  
The estimate consists of a term that is  a posteriori computable
from the solution, plus a remainder that is of higher order. 
It is the main tool for construction of   adaptive
algorithms. 
We present and test numerically one such algorithm.

The optimal control problem is  to minimize the functional
\begin{equation}
  \label{eq:objective}
  \int_0^T h(X(t),\alpha(t)) \ \ud t + g(X(T)),
\end{equation}
with given functions $h:\rset^{d}\times\mathcal{B}\to\rset$ and
$g:\rset^{d}\to\rset$, with respect to the state variable
$X:[0,T]\to\rset^{d}$ and the control
$\alpha:[0,T]\to\mathcal{B}$, with control set, $\mathcal{B}$, a
subset of some Euclidean space, $\rset^{d'}$, such that the ODE
constraint,
\begin{equation}
  \label{eq:constraint}
  \begin{aligned}
    X'(t) &= f(X(t),\alpha(t)), \quad 0 < t \leq T,\\
    X(0) &= x_0,    
  \end{aligned}
\end{equation}
is fulfilled.  This optimal control problem can be solved (globally)
using the Hamilton-Jacobi-Bellman (HJB) equation
\begin{equation}\label{eq:HJ}
  \begin{aligned}
    u_t + H(x,u_x) &= 0, & x\in\rset^{d}, \quad 0 \leq
    t < T,\\
    u(\cdot,T) &= g(\cdot), & x\in\rset^{d},\\
  \end{aligned}
\end{equation}
with $u_t$ and $u_x$ denoting the time derivative and spatial gradient
of $u$, respectively, and the Hamiltonian,
$H:\rset^{d}\times\rset^{d}\to\rset$, defined by
\begin{equation}\label{eq:Hamiltonian}
  H(x,\lambda) :=  \min_{\alpha\in\mathcal{B}} 
  \Bigl\{ 
  \lambda \cdot f(x,\alpha) +
  h(x,\alpha) 
  \Bigr\},
\end{equation}
and value function
\begin{equation}\label{eq:valuefunction}
  \begin{aligned}
    u(x,t) := \inf_{X:[t,T]\to\rset^{d},\ \alpha:[t,T]\to\mathcal{B}} \Biggl\{
    & \int_t^T h(X(s),\alpha(s)) \ \ud s +
    g(X(T)) \Biggr\},
  \end{aligned}
\end{equation}
where
\begin{equation*}
  \begin{aligned}
    X'(s) &= f(X(s),\alpha(s)), \quad t < s \leq T,\\
    X(t) &= x.
  \end{aligned}
\end{equation*}
The global minimum to the optimal control problem
\eqref{eq:objective}-\eqref{eq:constraint} is thus given by
$u(x_0,0)$.

If the Hamiltonian is sufficiently smooth, the bi-characteristics to the
HJB equation \eqref{eq:HJ} are given by the following Hamiltonian system:
\begin{equation}\label{eq:hamiltonian_system}
  \begin{aligned}
    X'(t) &= H_\lambda(X(t),\lambda(t)), &0 < t \leq T,\\
    X(0) &= x_0,\\
    -\lambda'(t) &= H_x(X(t),\lambda(t)), &0 \leq t < T,\\
    \lambda(T) &= g_x(X(T)),\\
  \end{aligned}
\end{equation}
where $H_\lambda$, $H_x$, and $g_x$ denote gradients with respect to
$\lambda$ and $x$, respectively, and the
dual variable, $\lambda:[0,T]\to\rset^{d}$, satisfies
$\lambda(t)= u_x(X(t),t)$ along the characteristic.

In Section \ref{sec:Error}, we present an error representation for
the following discretization of
\eqref{eq:hamiltonian_system}, which is used as a cornerstone for an
adaptive algorithm.
It is the Symplectic (forward) Euler method:
\begin{equation}\label{eq:symplectic}
  \begin{aligned}
    X_{n+1}- X_n &= \Delta t_n H_\lambda(X_n,\lambda_{n+1}), &
    n&=0,\ldots,N-1,\\
    X_0 &= x_0,\\
    \lambda_n-\lambda_{n+1} &= \Delta t_n 
    H_x(X_n,\lambda_{n+1}), & n&=0,\ldots,N-1,\\
    \lambda_N &= g_x(X_N),\\
  \end{aligned}
\end{equation}
with $0=t_0<t_1<\ldots<t_N=T$, $\Delta t_n:=t_{n+1}-t_n$, and
$X_n,\lambda_n\in\rset^{d}$. 
An alternative approach uses the dual weighted residual method, see
\cite{Becker-Rannacher,Bangerth-Rannacher}, to adaptively refine finite element solutions
of the Euler-Lagrange equation associated with the optimal control
problem, see \cite{KraftLarssonBIT,Kraft2011,Kraft2011a}.

The adaptive algorithm in Section \ref{sec:Error} uses a Hamiltonian
that is of $C^2$ regularity. In \cite{SS,Sandberg} first-order
convergence of the so-called Symplectic Pontryagin method, a
Symplectic Euler scheme \eqref{eq:symplectic} with a regularized
Hamiltonian $H^\delta$ replacing $H$, is shown. The Symplectic Pontryagin scheme
works in the more general optimal control setting where the
Hamiltonian is non-smooth. 
It uses the fact that if $u$ and $u^\delta$ are the solutions to the
Hamilton-Jacobi equation \eqref{eq:HJ} with the original (possibly
non-smooth) Hamiltonian $H$, and the regularized Hamiltonian,
$H^\delta$, then 
\begin{equation}\label{eq:udeltaerr}
  \norm{u-u^\delta}_{L^\infty([0,T]\times\rset^d)}\leq T
  \norm{H-H^\delta}_{L^\infty(\rset^d\times\rset^d)} 
  =\mathcal{O}(\delta),
\end{equation}
if $\norm{H-H^\delta}_{L^\infty(\rset^d\times\rset^d)}=\mathcal{O}(\delta)$.
Equation \eqref{eq:udeltaerr} is  a direct consequence of the maximum
principle for viscosity solutions to Hamilton-Jacobi equations, see
e.g., \cite{Bardi, Cannarsa-Sinestrari,Barles}.
For the error representation result in Theorem
\ref{thm:errorrepresentation}, we need $C^2$ regularity of
$H$. 
A possibility to use this error representation to find a solution
adaptively in the case where the Hamiltonian is non-differentiable, is to add
the error from the time discretization (the $\mathrm{TOL}$ in Theorem
\ref{thm:accuracy}) when the adaptive algorithm \ref{alg:adapt} is used
with a regularized Hamiltonian, $H^\delta$, to the error
$\mathcal{O}(\delta)$, in \eqref{eq:udeltaerr}. 
We show in Section \ref{sec:numerics}
that 
this method
works well for a test
case  in which the Hamiltonian is non-differentiable. 
Even though it works well in the cases we have studied, it is
difficult to justify this method theoretically. This is because the
size of the remainder term in Theorem \ref{thm:errorrepresentation}
depends on the size of the second-order derivatives of the Hamiltonian,
$H$, which typically are of order $\delta^{-1}$ when a regularized
$H^\delta$ is used.

\begin{rem}[Time-dependent Hamiltonian]\label{rem:timedep}
The analysis in this paper is presented for the optimal control
problem \eqref{eq:objective}, \eqref{eq:constraint}, i.e.,\ the case
where the running cost, $h$, and the flux, $f$, have no explicit time
dependence. The more general situation with explicit time dependence, 
to minimize 
  \begin{equation*}
    \int_0^T h(t,X(t),\alpha(t)) \ \ud t + g(X(T)),
  \end{equation*}
  for $\alpha\in\mathcal{B}$ such that the constraint
  \begin{equation*}
    \begin{aligned}
      X'(t) &= f(t,X(t),\alpha(t)), \quad 0 < t \leq T,\\
      X(0) &= x_0,
    \end{aligned}
  \end{equation*}
  is fulfilled, 
can be put in
  the form \eqref{eq:objective}, \eqref{eq:constraint} by introducing a
  state variable, $s(t)=t$, for the time dependence, i.e.,\ to minimize
    \begin{equation*}
    \int_0^T h(s(t),X(t),\alpha(t)) \ \ud t + g(X(T)),
  \end{equation*}
  such that the constraint
  \begin{equation*}
    \begin{aligned}
      X'(t) &= f(s(t),X(t),\alpha(t)), && 0 < t \leq T,\\
      s'(t) &= 1, && 0 < t \leq T,\\
      X(0) &= X_0,\\
      s(0) &= 0,
    \end{aligned}
  \end{equation*}
  is fulfilled. The Hamiltonian then becomes
  \begin{equation*}
    H(x,s,\lambda_1,\lambda_2) := \min_{\alpha\in\mathcal{B}} \Bigl\{
    \lambda_1\cdot f(x,\alpha,s) + \lambda_2 + h(x,\alpha,s) \Bigr\},
  \end{equation*}
where $\lambda_1$ is the dual variable corresponding to $X$, while
$\lambda_2$ corresponds to $s$.
\end{rem}

\section{Error estimation and adaptivity}\label{sec:Error}
In this section, we present an error representation for the Symplectic
Euler scheme in Theorem \ref{thm:errorrepresentation}. With this error
representation, it is possible to build an adaptive algorithm
(alg.\ \ref{alg:adapt}). The error representation in Theorem
\ref{thm:errorrepresentation} concerns approximation of the value
function, $u$, defined in \eqref{eq:valuefunction}. 
To define an approximate value function, $\bar u$, we need the
following definition of a running cost, a Legendre-type transform of
the Hamiltonian:
\begin{equation}\label{eq:LegendreL}
L(x,\beta)=\sup_{\lambda\in\mathbb{R}\sp d}
\big(-\beta\cdot\lambda+H(x,\lambda)\big),
\end{equation}
for all $x$ and $\beta$ in $\mathbb{R}\sp d$. The running cost
function is convex in its second argument and extended valued,
i.e.,\ its values belong to $\rset\cup\{+\infty\}$. 
If the Hamiltonian is real-valued and concave in its second variable,
it is possible to 
retrieve it from
$L$:
\begin{equation}\label{eq:LegendreH}
H(x,\lambda)=\inf_{\beta\in\mathbb{R}^d}\big(\lambda\cdot\beta+L(x,\beta)\big).
\end{equation}
This is a consequence of the bijectivity of the Legendre-Fenchel
transform, see \cite{Clarke, Sandberg}.

We now define a discrete value function:
\begin{equation}\label{eq:udef}
\bar u(y,t_m):=\inf\big\{J_{(y,t_m)}(\beta_m,\ldots,\beta_{N-1}) | \beta_m,\ldots,\beta_{N-1}\in\mathbb{R}^d\big\},
\end{equation}
where
\begin{equation}\label{eq:J}
J_{(y,t_m)}(\beta_m,\ldots,\beta_{N-1}) :=
\sum_{n=m}\sp{N-1} \Delta t_n L(X_n,\beta_n) +g(X_N),
\end{equation}
and 
\begin{equation}\label{eq:xalpharel}
\begin{split}
X_{n+1} &= X_n +\Delta t_n\beta_n,\quad\text{for }m \leq n\leq N-1,\\
X_m&=y.
\end{split}
\end{equation}
The appearance of a discrete path denoted $\{X_n\}$ in both the Symplectic
Euler scheme \eqref{eq:symplectic} and in the definition of $\bar u$
in \eqref{eq:xalpharel} is not just a coincidence.
The following theorem, taken from \cite{Sandberg}, shows that
to the minimizing path $\{X_n\}$ in the definition of $\bar u$
corresponds a discrete dual path $\{\lambda_n\}$, such that
$\{X_n,\lambda_n\}$ solves the Symplectic Euler scheme \eqref{eq:symplectic}.
For the statement and proof of Theorem \ref{thm:OptimEquivSP} we need
the following definitions.
\begin{defn}\label{def:semiconcavity} 
Let $S$ be a subset of $\Re\sp d$. We say that a function $f:S \rightarrow\Re$ is \emph{semiconcave} if
there exists a nondecreasing upper semicontinuous function $\omega:
\Re_+ \rightarrow \Re_+$ such that $\lim_{\rho\rightarrow 0\sp
  +}\omega(\rho) =0$ and
\begin{equation*}
w f(x)+(1-w) f(y) - f\big(w x+(1-w)y\big)\leq w(1-w)|x-y|\omega(|x-y|)
\end{equation*}
for any pair $x,y\in S$, such that the segment $[x,y]$ is contained in
$S$ and for any $w\in[0,1]$. We say that $f$ is \emph{locally
  semiconcave} on $S$ if it is semiconcave on every compact subset of $S$.
\end{defn}
There exist alternative definitions of semiconcavity, see
\cite{Cannarsa-Sinestrari}, but this is the one used in this paper.
\begin{defn}\label{def:superdifferential}
An element $p\in\Re\sp d$ belongs to the \emph{superdifferential} of
the function $f:\Re\sp d\rightarrow\Re$ at $x$, denoted $D\sp+ f(x)$, if 
\begin{equation*}
\limsup_{y\rightarrow x} \frac{f(y)-f(x)-p\cdot(y-x)}{|y-x|} \leq 0.
\end{equation*}
\end{defn}
\begin{thm}\label{thm:OptimEquivSP}
Let $y$ be any element in $\mathbb{R}\sp d$, and $g:\mathbb{R}\sp d\rightarrow\mathbb{R}$
a locally semiconcave function such that $g(x)\geq -k(1+|x|)$, for
some constant $k$, and all $x\in\mathbb{R}\sp d$. Let the Hamiltonian 
$H:\mathbb{R}\sp d\times\mathbb{R}\sp d\rightarrow \mathbb{R}$ satisfy the
following conditions:
\begin{itemize}
\item $H$ is differentiable everywhere in $\mathbb{R}\sp d\times\mathbb{R}\sp d$.
\item $H_\lambda(\cdot,\lambda)$ is locally Lipschitz continuous for every $\lambda \in \mathbb{R}\sp d$.
\item $H_x$ is continuous everywhere in $\mathbb{R}\sp d\times\mathbb{R}\sp d$.
\item There exists a convex, nondecreasing function
  $\mu:[0,\infty)\rightarrow \mathbb{R}$ and positive constants $A$ and $B$
  such that 
\begin{equation}\label{eq:Hmu}
  -H(x,\lambda)\leq \mu(|\lambda|) + |x|(A+B|\lambda|)\quad
   \text{for all }  (x,\lambda)\in\mathbb{R}\sp d\times\mathbb{R}\sp d.
\end{equation}
\item $H(x,\cdot)$ is concave for every $x \in \mathbb{R}\sp d$.
\end{itemize}
Let 
$L$ be defined
by
\eqref{eq:LegendreL}. Then, there exists a minimizer
$(\beta_m,\ldots,\beta_{N-1})$ of the function $J_{(y,t_m)}$ in \eqref{eq:J}. Let
$(X_m,\ldots,X_N)$ be the corresponding solution to \eqref{eq:xalpharel}.
Then, for each $\lambda_N\in D\sp + g(X_N)$, there exists a discrete dual path
$(\lambda_m,\ldots,\lambda_{N-1})$, that satisfies 
\begin{equation}\label{eq:SymplPontr}
\begin{split}
X_{n+1}&=X_n+\Delta t_n H_\lambda(X_n,\lambda_{n+1}), \quad\text{for all
} m \leq n\leq N-1,\\
X_m&=y\\
\lambda_n &= \lambda_{n+1}+\Delta t_n
H_x(X_n,\lambda_{n+1}),\quad\text{for all }m \leq n\leq N-1.
\end{split}
\end{equation}
Hence, 
\begin{equation}\label{eq:betan}
\beta_n=H_\lambda(X_n,\lambda_{n+1})
\end{equation}
for all $m \leq n\leq N-1$.
\end{thm}

The proof of Theorem \ref{thm:OptimEquivSP} from \cite{Sandberg} is
reproduced in the appendix.
  
With the correspondence between the Symplectic Euler scheme and discrete
minimization in Theorem \ref{thm:OptimEquivSP}, we are now ready to
formulate the error representation result. We will use the terminology
that a function is bounded in $C^k$ if it belongs to $C^k$ and has
bounded derivatives of order less than or equal to $k$. 
\begin{thm}\label{thm:errorrepresentation}
Assume that all conditions in Theorem \ref{thm:OptimEquivSP}
are satisfied,
that the Hamiltonian, $H$, is bounded
in $C^2(\mathbb{R}^d\times\mathbb{R}^d)$,
and that there exists a constant, $C$, such that for every discretization $\{t_n\}$ the difference between
the discrete dual and the gradient of the value function is bounded as
\begin{equation*}
\abs{\lambda_n-u_x(X_n,t_n)}\leq C\Delta t_{\mathrm{max}},
\end{equation*}
where $\Delta t_{\mathrm{max}}:=\max_n \Delta t_{n}$.
Assume further that either of the following two conditions holds:
\begin{enumerate}
\item\label{cond:1} The value function, $u$, is bounded in $C^3((0,T)\times\mathbb{R}^d)$. 
\item\label{cond:2} There exists a neighborhood in $C([0,T],\rset^d)$ around the minimizer
  $X:[0,T]\rightarrow\rset^d$ of $u(x_0,0)$ in
  \eqref{eq:valuefunction} in which the value function, $u$, is
  bounded in $C^3$. Moreover, the discrete solutions $\{X_n\}$
  converge to the continuous solution $X(t)$ in the sense that
\begin{equation*}
\max_n\abs{X_n-X(t_n)}\rightarrow 0,\ \text{as } \Delta
t_{\mathrm{max}}\rightarrow 0.
\end{equation*}
\end{enumerate}

If Condition \ref{cond:1} holds, then for every discretization $\{t_n\}$,
the error $\bar u(x_0,0)-u(x_0,0)$ is given as 
\begin{equation}\label{eq:errorrepr}
  \bar u(x_0,0)-u(x_0,0)=\sum_{n=0}^{N-1} \Delta
  t_n^2\rho_n  +  R,
\end{equation}
with density
\begin{equation}\label{eq:density}
\rho_n:=-\frac{H_\lambda(X_n,\lambda_{n+1})\cdot H_x(X_n,\lambda_{n+1})}{2}
\end{equation}
and 
the remainder term, $\abs{R}\leq C'\Delta t_{\mathrm{max}}^2$, for some
constant $C'$.

If Condition \ref{cond:2} holds, then there exists a threshold time
step, $\Delta t_{\mathrm{thres}}$, such that for every discretization
with $\Delta t_{\mathrm{max}}\leq \Delta t_{\mathrm{thres}}$ the error
representation \eqref{eq:errorrepr} holds. 
\end{thm}
\begin{rem}
In the proof of the 
theorem, we show that equation \eqref{eq:errorrepr} is satisfied with   the error density 
\begin{equation}\label{eq:tildedensity}
  \begin{aligned}
    \tilde\rho_n :=&\ \frac{H(X_n,\lambda_{n + 1})}{\Delta
      t_n}-\frac{H(X_n,\lambda_n) + H(X_{n + 1},\lambda_{n + 1})}{2\Delta
      t_n}  \\
    &+\frac{\lambda_n-\lambda_{n + 1}}{2}\cdot
    \frac{H_\lambda(X_n,\lambda_{n + 1})}{\Delta t_n}
  \end{aligned}
\end{equation}
replacing $\rho_n$. 
Under the assumption that the Hamiltonian, $H$, is bounded in $C^2$, we
have that $\abs{\rho_n-\tilde\rho_n}=\mathcal{O}(\Delta t_n)$. This
follows by Taylor expansion and by using that $\{X_n,\lambda_n\}$ solves
the Symplectic Euler scheme \eqref{eq:symplectic}. Hence, the theorem holds also with the error density
$\rho_n$. An advantage of $\rho_n$ is that 
it is given by a simple
expression. The error density $\tilde\rho_n$ has the advantage that it
is the one that is obtained in the proof, and then $\rho_n$ is derived
from it. One could therefore expect that $\tilde\rho_n$ would give a
slightly more accurate error representation. Moreover, $\tilde\rho_n$
is directly computable (as is $\rho_n$) once a solution
$\{X_n,\lambda_n\}$ has been computed.
\end{rem}
\begin{proof}
By Theorem \ref{thm:OptimEquivSP},
the error can be expressed as 
\begin{equation}\label{eq:err1}
  (\bar u-u)(x_0,0) = \sum_{n=0}^{N-1} \Delta t_n L(X_n,\beta_n) +
  g(X_N) - u(x_0,0),
\end{equation}
where 
\begin{equation*}
  g(X_N)=u(X_N,T),\quad\beta_n = H_\lambda(X_n,\lambda_{n+1}).
\end{equation*}

Define the piecewise linear function $\bar X(t)$ to be
\begin{equation*}
  \begin{aligned}
    \bar X(t) &= X_n + (t-t_n)H_\lambda(X_n,\lambda_{n+1}), \qquad
    t\in(t_n,t_{n+1}), \quad n=0,\ldots,N-1.
  \end{aligned}
\end{equation*}
If Condition \ref{cond:2} in the theorem holds, we now assume that
$\Delta t_{\mathrm{max}}$ is small enough, such that the path $\bar X(t)$
belongs to the neighborhood of $X(t)$ in $C([0,T],\rset^d)$ where the
value function belongs to $C^3$.
If Condition \ref{cond:1} holds, the following analysis is also valid,
without restriction on $\Delta t_{\mathrm{max}}$.
From \eqref{eq:err1} and the Hamilton-Jacobi-Bellman equation, we have
\begin{equation}\label{eq:err2}
  \begin{aligned}
    (\bar u-u)(x_0,0) =&\ \sum_{n=0}^{N-1} \Delta t_n L(X_n,\beta_n) +
    u(X_N,T) -
    u(x_0,0)\\
    =&\ \sum_{n=0}^{N-1} \Delta t_n L(X_n,\beta_n) + \int_0^T
    \frac{\ud }{\ud
      t}u(\bar X(t),t)\, \ud t\\
    =&\ \sum_{n=0}^{N-1}  \int_{t_n}^{t_{n+1}}  L(X_n,\beta_n)\, \ud t \\
    &+ \sum_{n=0}^{N-1} \int_{t_n}^{t_{n+1}} u_t(\bar
      X(t),t) +
    u_x(\bar X(t),t)\cdot
    H_\lambda(X_n,\lambda_{n+1})\, \ud t.
  \end{aligned}
\end{equation}
By \eqref{eq:LegendreL} and \eqref{eq:betan} we have
\begin{equation*}
  H(X_n,\lambda_{n+1}) = \lambda_{n+1} \cdot
  H_\lambda(X_n,\lambda_{n+1}) + L(X_n,\beta_n),
\end{equation*}
which 
together with the Hamilton-Jacobi equation
\begin{equation*}
u_t(\bar X(t),t)=-H\big(\bar X(t),u_x(\bar X(t),t)\big)
\end{equation*}
implies that the error can be written as
\begin{equation}\label{eq:Errorsum}
  \begin{aligned}
    (\bar u-u)(x_0,0) &= \sum_{n=0}^{N-1} \int_{t_n}^{t_{n + 1}}
    H(X_n,\lambda_{n + 1})-H(\bar
    X(t), u_x(\bar X(t),t)) \, \ud t  \\
    &\quad +\sum_{n=0}^{N-1} \int_{t_n}^{t_{n + 1}} \big(u_x(\bar
    X(t),t)-\lambda_{n + 1}\big)\cdot H_\lambda(X_n, \lambda_{n + 1})
    \, \ud t\\
    &=: \sum_{n=0}^{N-1} E_n.
  \end{aligned}
\end{equation}
By the boundedness of the Hamiltonian, $H$, in $C^2$ and the
value function, $u$, in $C^3$, it follows that the trapezoidal
rule can be applied to the integrals in \eqref{eq:Errorsum} with an
error of order $\Delta t_n^3$.  Hence, we obtain that
\begin{equation}\label{eq:Errorsum2}
  \begin{aligned}
    E_n =&\ \Delta t_n\Bigl( H(X_n,\lambda_{n + 1})-\frac{H(X_n,
      u_x(X_n,t_n)) + H(X_{n + 1},u_x(X_{n + 1},t_{n + 1}))}{2} \Bigr)\\
    &+ \Delta t_n \Bigl(\frac{u_x(X_n,t_n) + u_x(X_{n
        + 1},t_{n + 1})}{2}-\lambda_{n + 1}\Bigr)\cdot H_\lambda(X_n,
    \lambda_{n + 1}) + \bar R_n
  \end{aligned}
\end{equation}
with  remainder $\bar R_n=\mathcal{O}(\Delta
t_{n}^3)$. 

What remains for us to show is that we can exchange the gradient of the
continuous value function, $u$, in \eqref{eq:Errorsum2} with the
discrete dual, $\lambda_n$, with an error bounded by $\Delta
t_{\mathrm{max}}^2$. 
We write this difference using the 
error density, $\tilde\rho_n$, from \eqref{eq:tildedensity}: 
\begin{equation*}\label{eq:difference}
  \begin{aligned}
    \Delta t_n^2 \tilde\rho_n - &E_n =\  - \frac{\Delta t_n}{2}
    \Bigl( H(X_n,\lambda_{n})-H(X_n,u_x(X_n,t_n))\Bigr) \\
    &- \frac{\Delta t_n}{2} \Bigl( H(X_{n + 1},\lambda_{n + 1})-H(X_{n
      + 1},u_x(X_{n +
      1},t_{n + 1})) \Bigr) \\
    &+ \frac{\Delta t_n }{2} \Bigl( \lambda_n - u_x(X_n,t_n) +
    \lambda_{n + 1} - u_x(X_{n + 1},t_{n + 1}) \Bigr) \cdot
    H_\lambda(X_n,\lambda_{n + 1}) - \bar R_n\\
    =&\  \frac{\Delta t_n}{2} \Bigl( -E^I_n - E^{I}_{n+1} +
    (\xi_n+\xi_{n+1})\cdot H_\lambda(X_n,\lambda_{n
      + 1} \Bigr) - \bar R_n ,
  \end{aligned}
\end{equation*}
where
\begin{gather*}
    E^I_n :=\  H(X_n,\lambda_{n})-H(X_n,u_x(X_n,t_n))
    =\ H_\lambda(X_n,\lambda_{n}) \cdot \xi_n +
    \mathcal{O}\Bigl( |\xi_n|^2 \Bigr),\\
  \xi_n :=  \lambda_n- u_x(X_n,t_n).
\end{gather*}
Further Taylor expansion gives the difference
\begin{equation*}
  \begin{aligned}
    E^I_n - \xi_n \cdot H_\lambda(X_n,\lambda_{n + 1}) =&\
    \Bigl( H_\lambda(X_n,\lambda_{n}) - H_\lambda(X_n,\lambda_{n +
      1})\Bigr) \cdot \xi_n +
    \mathcal{O}\Bigl( |\xi_n|^2 \Bigr)\\
    =&\ \mathcal{O}\Bigl( \Delta t_n |\xi_n| +
    |\xi_n|^2\Bigr) = \mathcal{O}\Bigl( \Delta
    t_{\mathrm{max}}^2\Bigr),
  \end{aligned}
\end{equation*}
and similarly
\begin{equation*}
  \begin{aligned}
    E^I_{n+1} - \xi_{n+1} \cdot H_\lambda(X_n,\lambda_{n +
      1}) = \mathcal{O}\Bigl( \Delta t_{\mathrm{max}}^2\Bigr).
  \end{aligned}
\end{equation*}

Finally, summing the difference $\Delta t_n^2 \tilde\rho_n - E_n$ over
$n=0,\ldots,N-1$ gives, together with the above Taylor expansions, the
bound $|R|\leq C \Delta t_{\mathrm{max}}^2$ in the theorem.

\end{proof}

In what follows, we formulate an adaptive algorithm (\ref{alg:adapt}) and
three theorems (\ref{thm:stopping}--\ref{thm:efficiency}) on its
performance. These are all taken from \cite{Moon2003c} more or less
directly. Since the proofs are practically unchanged, they are not
repeated here.

\begin{alg}[Adaptivity]\label{alg:adapt}
  Choose the error tolerance $\mathrm{TOL}$, the 
  initial grid
  $\{t_n\}_{n=0}^N$, the parameters $s$ and $M$, 
  and repeat
  the following points:
  \begin{enumerate}
  \item Calculate $\{(X_n,\lambda_n)\}_{n=0}^N$ with the symplectic Euler scheme \eqref{eq:symplectic}.
  \item Calculate error densities $\{\rho_n\}_{n=0}^{N-1}$ and the
    corresponding approximate error densities
    \begin{equation*}
      \bar \rho_n := \sgn(\rho_n)\max(\abs{\rho_n},K\sqrt{\Delta t_{\mathrm{max}}}).
    \end{equation*}
  \item Break if
    \begin{equation*}
      \max_n \bar r_n < 
      \frac{\mathrm{TOL}}{N}
    \end{equation*}
    where the error indicators are defined by $\bar r_n:=\abs{\bar
    \rho_n}\Delta t_n^2$.
  \item Traverse through the mesh and subdivide an interval
    $(t_n,t_{n+1})$ into $M$ parts if
    \begin{equation*}
      \bar r_n > s \frac{\mathrm{TOL}}{N}.
    \end{equation*}
  \item Update $N$ and $\{t_n\}_{n=0}^N$ to reflect the new mesh.
  \end{enumerate}

\end{alg}
The goal of the algorithm is to construct a partition of the time
interval $[0,T]$ such that 
\begin{equation*}
\bar r_n\approx \frac{\mathrm{TOL}}{N},
\end{equation*}
for all $n$. The constant $s<1$ is present in order to achieve a
substantial reduction of the error, described further in Theorem \ref{thm:stopping}. 
The constant $K$ in the algorithm should be chosen small (relative to
the size of the solution). In the numerical experiments presented in
Section \ref{sec:numerics}, we use $K=10^{-6}$.

Let $\Delta t(t)[k]$ be defined as the piecewise constant function that
equals the local time step 
\begin{equation*}
\Delta t(t)=\Delta t_n,\quad \text{if } t\in[t_n,t_{n+1}),
\end{equation*}
on mesh refinement level $k$.
As in \cite{Moon2003c}, we have that 
\begin{equation*}
\lim_{\mathrm{TOL}\rightarrow0^+}\max_t \Delta t(t)[J]=0,
\end{equation*}
where mesh $J$ is the finest mesh where the algorithm stops. By the
assumptions on the convergence of the approximate paths
$\{X_n,\lambda_n\}$, it follows that there exists a limit
\begin{equation*}
\abs{\bar\rho}\rightarrow\abs{\tilde\rho},\ \text{ as } \max\Delta
t\rightarrow 0.
\end{equation*}
We introduce a constant, $c=c(t)$, such that
\begin{equation}\label{eq:rhobounds} 
\begin{aligned}
c&\leq\abs{\frac{\bar\rho(t)[\text{parent}(n,k)]}{\bar\rho(t)[k]}}\leq
  c^{-1},\\ 
c&\leq\abs{\frac{\bar\rho(t)[k-1]}{\bar\rho(t)[k]}}\leq
  c^{-1},\\ 
\end{aligned}
\end{equation}
holds for all time steps, $t\in\Delta t_n[k]$, and all refinement levels,
  $k$. Here, $\text{parent}(n,k)$ means the refinement level where a
coarser interval was split into a number of finer subintervals of
which $\Delta t_n[k]$ is one.
Since $\abs{\bar\rho}$ converges as $\mathrm{TOL}\rightarrow 0$ and is
  bounded away from zero, $c$ will be close to $1$ for sufficiently
  fine meshes.

\begin{thm}\label{thm:stopping}[Stopping]
 Assume that $c$ satisfies \eqref{eq:rhobounds} for the time steps
 corresponding to the maximal error indicator on each refinement
 level, and that 
\begin{equation}\label{eq:Stoppingparameters}
M^2>c^{-1},\ 
s\leq\frac{c}{M}.
\end{equation}
Then, each refinement level either decreases the maximal error
indicator with the factor
\begin{equation*}
\max_{n}\bar
r_n[k+1]\leq\frac{c^{-1}}{M^{2}}\max_{n}\bar r_n[k],
\end{equation*}
or stops the algorithm.
\end{thm}

The inequalities in \eqref{eq:Stoppingparameters} give (at least in
principle) an idea how to determine the parameters $M$ and $s$. When
the constant, $c=c(t)$, has been determined approximately, say after one
or a few refinements, $M$  can be chosen using the first inequality
and then $s$ can be chosen using the other. 
\begin{thm}\label{thm:accuracy}[Accuracy]
  The adaptive Algorithm \ref{alg:adapt} satisfies
\begin{equation*}
\limsup_{\mathrm{TOL}\rightarrow 0^+}\bigl(\mathrm{TOL}^{-1}\abs{u(x_0,0)-\bar
  u(x_0,0)}\bigr)\leq 1.
\end{equation*} 
\end{thm}

\begin{thm}\label{thm:efficiency}[Efficiency]
Assume that $c=c(t)$ satisfies \eqref{eq:rhobounds} for all time steps
at the final refinement level, and that all initial time steps have been
divided when the algorithm stops. Then, there exists a constant, $C>0$,
bounded by $M^2s^{-1}$, such that the final number of adaptive
steps, $N$, of the Algorithm \ref{alg:adapt}, satisfies
\begin{equation*}
\mathrm{TOL}\ N\leq
C\norm{\frac{\bar\rho}{c}}_{L^{\frac{1}{2}}}\leq
\norm{\bar\rho}_{L^{\frac{1}{2}}} \max_{0\leq
  t\leq T} c(t)^{-1},
\end{equation*}
and $\norm{\bar\rho}_{L^{\frac{1}{2}}}\rightarrow
\norm{\tilde\rho}_{L^{\frac{1}{2}}}$ asymptotically as $\mathrm{TOL}\rightarrow
0^+$.
\end{thm}
\begin{rem}
Note that the optimal number $N_a$ of non-constant (i.e., adaptive) time
steps to have the error $\sum_n \Delta t_n^2\bar\rho_n$ smaller than
$\mathrm{TOL}$ satisfies $\mathrm{TOL} N_a\approx \|\bar\rho\|_{L^{1/2}}$, see
\cite{Moon2003c}, while the number of uniform time steps $N_u$ required satisfies $\mathrm{TOL} N_u\approx \|\bar\rho\|_{L^{1}}$.
\end{rem}
\begin{rem}
It is natural to use adaptivity when  optimal control problems are
solved using the Hamiltonian system
\eqref{eq:hamiltonian_system}. Since it is a coupled ODE system with a
terminal condition linking the primal and dual functions, it is
necessary to solve using some iterative method. When an initial guess
is to be provided to the iterative method, it is natural to interpolate
a solution obtained on a coarser mesh. Solutions on several
meshes therefore need be computed, as is the case when adaptivity is used.
\end{rem}

\section{Numerical examples}\label{sec:numerics}
In this section, we consider three numerical examples. The first
is an optimal control problem that satisfies the assumption of a
$C^2$ Hamiltonian in Theorem \ref{thm:errorrepresentation}. The second
is a problem in which the Hamiltonian is non-differentiable, and hence
does not fulfill the smoothness assumption of Theorem
\ref{thm:errorrepresentation}. We investigate the influence of a
regularization of the Hamiltonian. The third example is a problem
in which the controlled ODE  has an explicit time dependence with a
singularity.

We will compare the work and error for the
adaptive mesh refinement in Algorithm \ref{alg:adapt} with that of
uniform mesh refinement. The work is represented by the cumulative
number of time steps on all refinement levels,
and the error is represented by either an estimation of the true
error, using the value function from the finest unform mesh as our
true solution, or estimating the error by
\begin{equation}\label{eq:error}
  E := \abs{\sum_{n=1}^{N-1} \bar \rho_n \Delta t_n^2},
\end{equation}
using the approximate error densities,
\begin{equation*}
  \bar \rho_n := \sgn(\rho_n)\max(\abs{\rho_n},10^{-6}\sqrt{\Delta t_{\mathrm{max}}}).
\end{equation*}

In all examples, we let $s=0.25$ and $M=2$ (since $c\approx 1$). 
On each mesh, the discretized Hamiltonian system \eqref{eq:symplectic}
is solved with MATLAB's {\tt FSOLVE} routine, with default parameters
and a user-supplied Jacobian, and using the solution from the previous
mesh as a starting guess.

\begin{ex}[Hyper-sensitive optimal control]\label{ex:hyper}
  This is a version of Example 6.1 in \cite{Kraft2011a} and Example 51
  in \cite{propt}. Minimize
  \begin{equation*}
    \int_0^{25} \big(X(t)^2+\alpha(t)^2\big)\, \ud t + \gamma (X(25)-1)^2,
  \end{equation*}
  subject to
  \begin{equation*}
    \begin{aligned}
      X'(t) &= -X(t)^3+\alpha(t), \qquad 0<t\leq 25,\\
      X(0) &= 1,
    \end{aligned}
  \end{equation*}
  for some large $\gamma>0$.
  The Hamiltonian is then given by
  \begin{equation*}
    \begin{aligned}
      H(x,\lambda) :=&\ \min_\alpha \Bigl\{ -\lambda x^3 + \lambda
      \alpha + x^2 + \alpha^2 \Bigr\} = -\lambda x^3 - \lambda^2/4 +
      x^2.
    \end{aligned}
  \end{equation*}
\end{ex}

First,
we run the adaptive algorithm with tolerance, $\mathrm{TOL}$, leading to the
estimated error, $E_{\mathrm{adap}}$. Finally, the problem is rerun using uniform
refinement with stopping criteria, $E_{\mathrm{unif}}\leq E_{\mathrm{adap}}$.

\begin{figure}
  \includegraphics[width=0.99\textwidth]
  {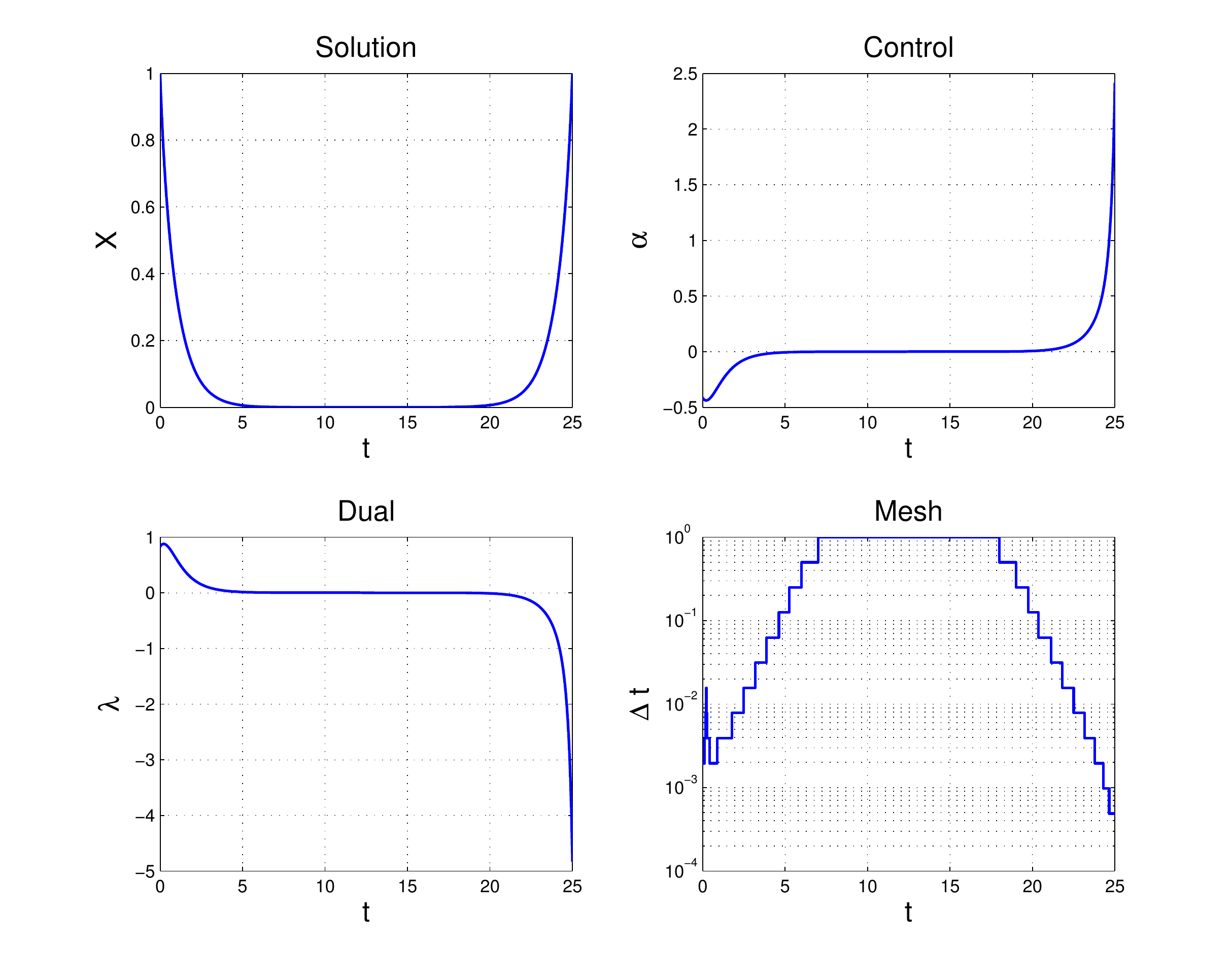}
  \caption{The solution, $X$, control, $\alpha$, dual, $\lambda$, and mesh,
    $\Delta t$, for the hyper-sensitive optimal control problem in Example
    \ref{ex:hyper}, with $\gamma=10^{6}$ and $\mathrm{TOL}=10^{-2}$.
   }  
 \label{fig:hyper1}
\end{figure}

\begin{figure}
  \includegraphics[width=0.99\textwidth]
  {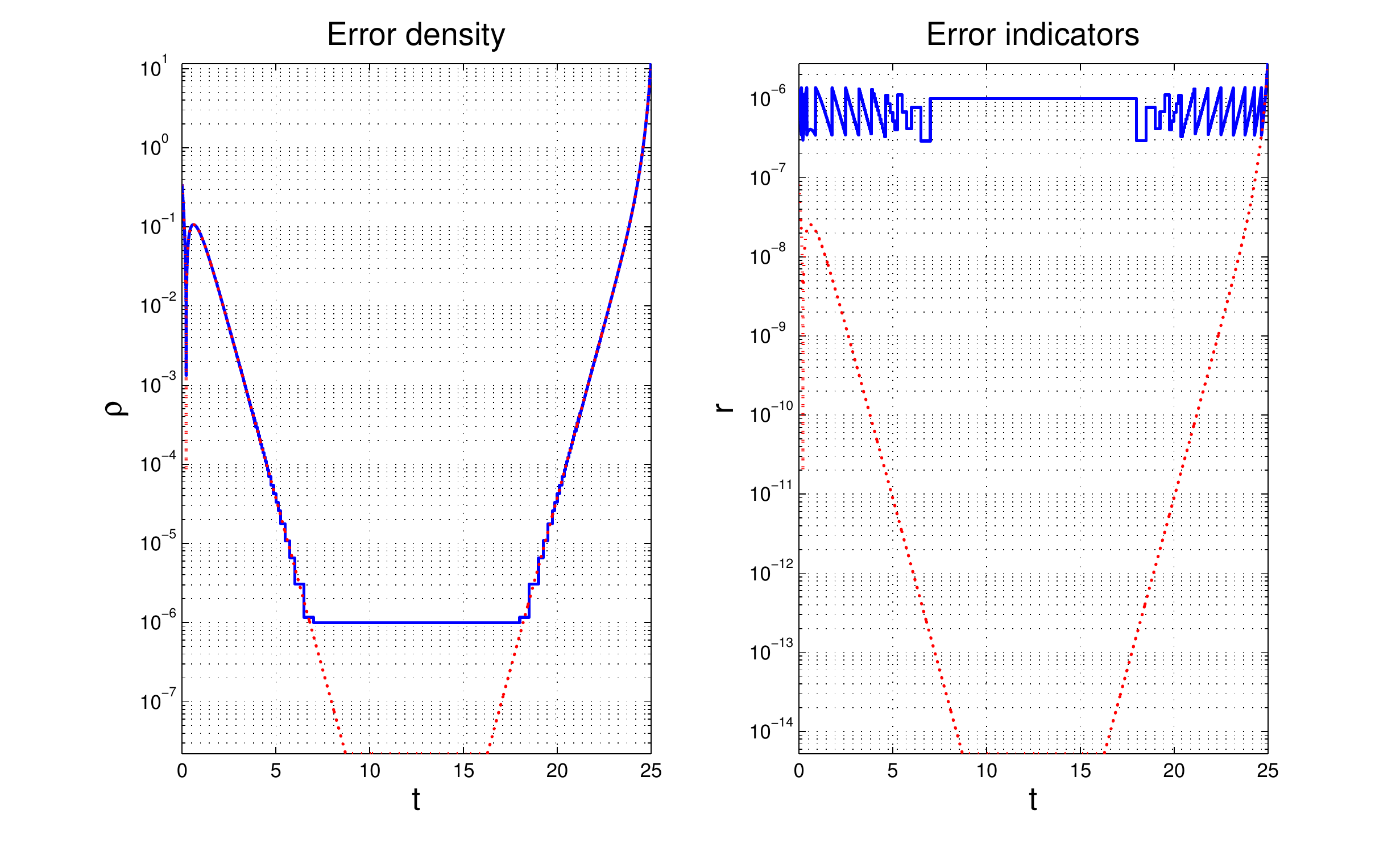}
  \caption{Error densities, $\abs{\bar \rho_n}$, and error indicators, $\bar
    r_n$, for the hyper-sensitive optimal control problem in Example
    \ref{ex:hyper}. The solid and dotted lines correspond to solutions with adaptive and uniform time stepping, respectively.}
  \label{fig:hyper2}
\end{figure}

\begin{figure}
\centering
  \includegraphics[width=0.7\textwidth]
  {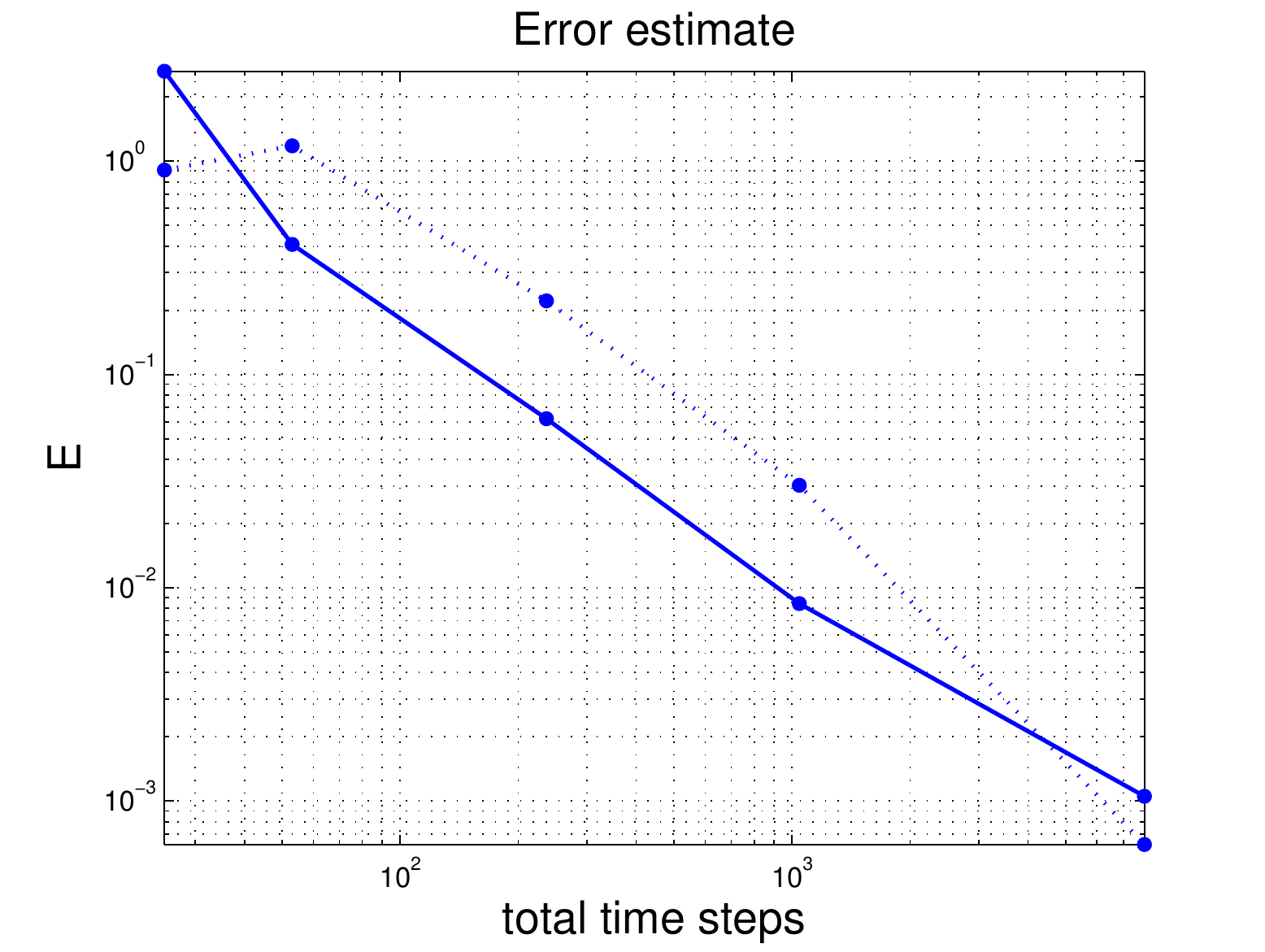}
  \caption{Error estimates for the hyper-sensitive optimal control
    problem in Example \ref{ex:hyper}. The solid line indicates the
    error estimate in \eqref{eq:error}, and the dotted
    line indicates the difference between the
    value function and the value function using a fine uniform mesh
    with 51200 time steps. The error  estimate from 
    \eqref{eq:error} for the uniform mesh is approximately as large as
  the estimate for the finest adaptive level. Hence, the dotted line is
  only an approximation of the true error.}
  \label{fig:hypererrorcomparison}
\end{figure}

\begin{figure}
\centering
  \includegraphics[width=0.7\textwidth]
  {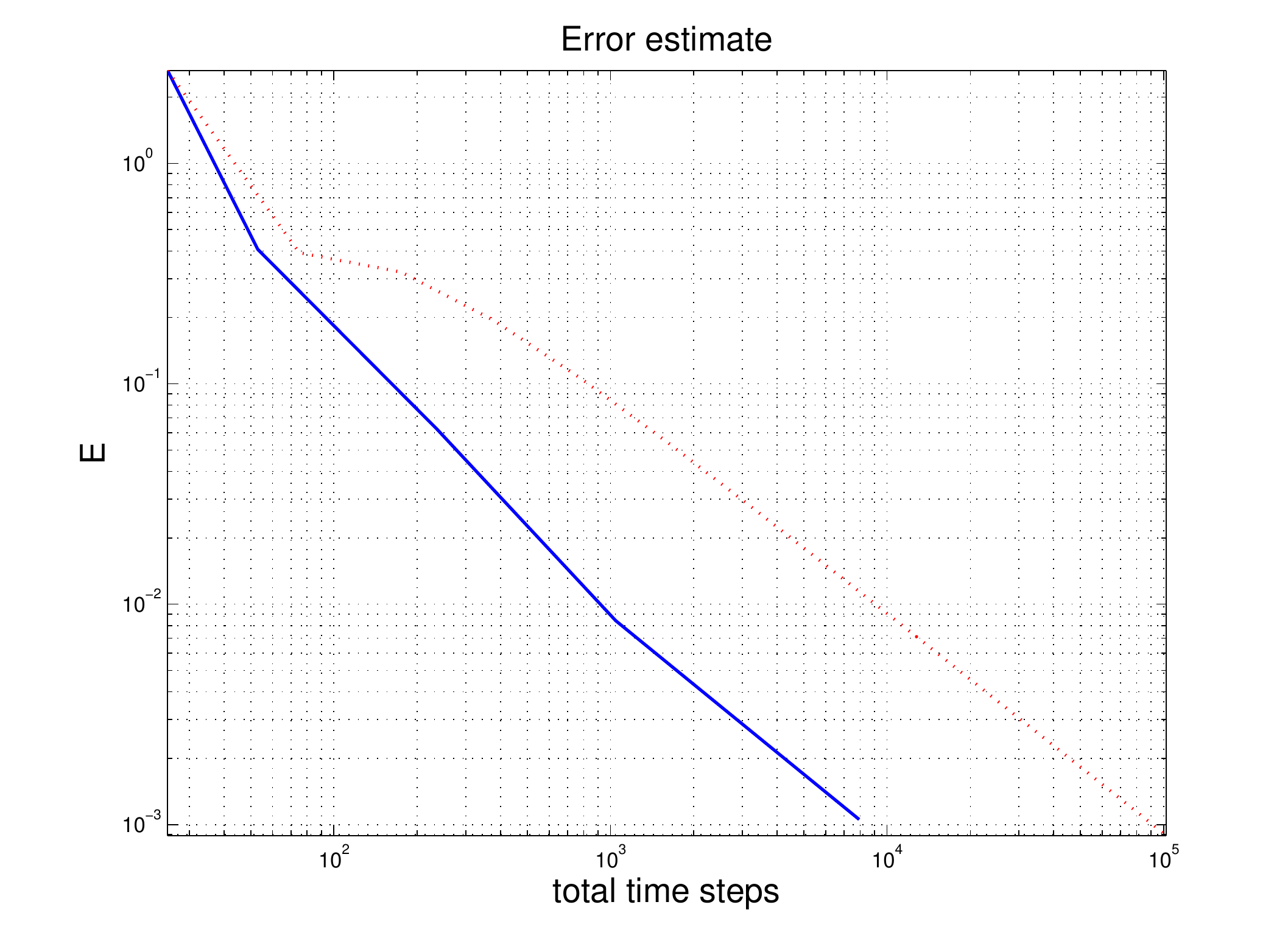}
  \caption{Error estimates for the hyper-sensitive optimal control
    problem in Example \ref{ex:hyper} using \eqref{eq:error}, versus the 
    cumulative number of time steps on all refinement levels for the adaptive algorithm (solid) and
    uniform meshes (dotted). The number of time steps in the uniform
    meshes is doubled in each refinement.} 
  \label{fig:hyperworkcomparison}
\end{figure}

Figure \ref{fig:hyper1} shows the solution and final mesh when
computed with the adaptive
Algorithm \ref{alg:adapt}.
Figure \ref{fig:hyper2} shows the error density and error indicator,
while Figure \ref{fig:hypererrorcomparison} gives a comparison between
the error estimate from equation \eqref{eq:error} with an estimate of
the error using a uniform mesh solution with a small step size as a reference. 
Figure \ref{fig:hyperworkcomparison} shows error estimates versus
computational work as the cumulative number of time steps.

The error representation in Theorem \ref{thm:errorrepresentation}
concerns approximation of the value function when  the Symplectic
Euler scheme is used  with a $C^2$ Hamiltonian. In general, the
minimizing $\alpha$ in the definition of the Hamiltonian
\eqref{eq:Hamiltonian} depends discontinuously on $x$ and $\lambda$,
which most probably leads to a non-differentiable Hamiltonian. In
Example \ref{ex:simple} 
we consider a simple optimal control problem with an
associated Hamiltonian that is  non-differentiable. We denote by
$H^\delta$ a $C^2$ regularization of the Hamiltonian, $H$, such that
\begin{equation*}
\norm{H-H^{\delta}}_{L^{\infty}(\rset^d\times\rset^d)}=\mathcal{O}(\delta).
\end{equation*}
Since the remainder term in Theorem \ref{thm:errorrepresentation}
contains second-order derivatives of the Hamiltonian, which are of order
$\delta^{-1}$ if a regularization $H^\delta$ is used, it could be
expected that an estimate of the error using the error density term 
\begin{equation}\label{eq:approximateerror}
\sum_{n=0}^{N-1}\Delta t_n^2\rho_n
\end{equation}
in \eqref{eq:errorrepr} would be imprecise. However, the solution of
Example \ref{ex:simple} 
suggests that the approximation of the error in
\eqref{eq:approximateerror} might be accurate even in cases where
regularization is needed and the regularization parameter, $\delta$, is
chosen to be small.

\begin{ex}[A simple optimal control problem]\label{ex:simple}
Minimize
  \begin{equation}\label{eq:simplefunctional}
    \int_0^1 X(t)^{10}\, \ud t,
  \end{equation}
  subject to
  \begin{equation*}
    \begin{aligned}
      X'(t) &= \alpha(t) \in [-1,1], & 0<t\leq T,\\
      X(0) &= 0.5.\\
    \end{aligned}
  \end{equation*}
  The Hamiltonian is then non-smooth:
  \begin{equation*}
    \begin{aligned}
      H(x,\lambda) :=&\ \min_{\alpha\in[-1,1]} \Bigl\{ \lambda \alpha
      + x^{10} \Bigr\} = -\abs{\lambda} + x^{10},
    \end{aligned}
  \end{equation*}
  but can be regularized by
  \begin{equation*}
    \begin{aligned}
      H_\delta(x,\lambda) :=&\ - \sqrt{\lambda^2+\delta^2} + x^{10},
    \end{aligned}
  \end{equation*}
  for some small $\delta>0$.
\end{ex}

  The exact solution, without regularization, is $X(t)=(0.5-t)$ for
  $t\in[0,0.5]$ and $X(t)=0$ elsewhere, with control $\alpha(t)=-1$
  for $t\in[0,0.5]$ and $\alpha(t)=0$ elsewhere. This gives the
  optimal  value of the cost functional \eqref{eq:simplefunctional} (the value function) to be 
  $0.5^{11}/11$. 

\begin{figure}
\centering
  \includegraphics[width=0.7\textwidth]
  {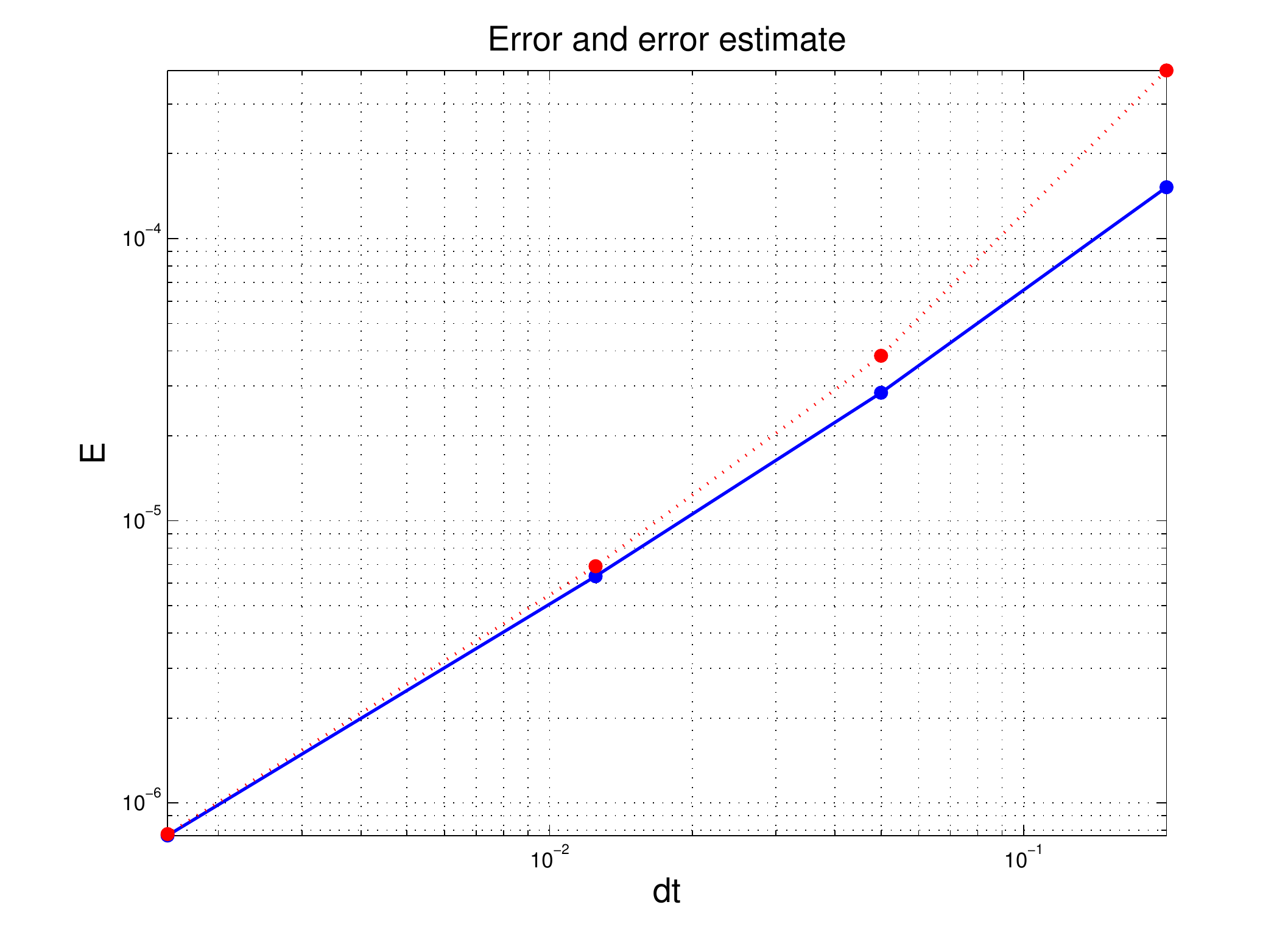}
  \caption{The true error (solid) and error estimation using
    \eqref{eq:approximateerror} (dotted) for the simple optimal control problem in Example \ref{ex:simple}
    with regularization parameter $\delta=10^{-10}$.}
  \label{fig:simplecomparison}
\end{figure}

  In Figure \ref{fig:simplecomparison}, a comparison is made between
  the error estimate, $\sum_{n=0}^{N-1}\Delta t_n^2\rho_n$, and the true
  error. It seems clear that the error estimate converges to the true
  error as $\Delta t\rightarrow 0$. In this numerical test, the
  regularization parameter, $\delta=10^{-10}$, and hence the part of
  the error from the regularization is negligible.
\begin{ex}[A singular optimal control problem]\label{ex:singular}
  This example is based on the singular ODE example in
  \cite{Moon2003c}, suitable for adaptive refinement. Consider the optimal control problem to minimize
  \begin{equation}\label{eq:singularminimization}
    \int_0^4\big(\alpha(t)-X(t)\big)^2 \ud t + \big(X(4)-X_\mathrm{ref}(4)\big)^2
  \end{equation}
  under the constraint
  \[
  \begin{split}
  X'(t)&=\frac{\alpha(t)}{\big((t-t_0)^2+\varepsilon^2\big)^{\beta/2}},\\
  X(0)&=X_\mathrm{ref}(0),
  \end{split}
  \]
  where $t_0=5/3$. 
  The reference $X_{\mathrm{ref}}(t)$ solves 
  \[
  X_\mathrm{ref}'(t)=\frac{X_\mathrm{ref}(t)}{\big((t-5/3)^2+\varepsilon^2\big)^{\beta/2}}
  \]
  and is given explicitly by
  \[
  X_{\mathrm{ref}}(t)=\exp\Big(\frac{t-t_0}{\varepsilon^\beta}\, {}_2F_1(\frac{1}{2},\frac{\beta}{2},\frac{3}{2};-\frac{(t-t_0)^2}{\varepsilon^2})\Big),
  \]
  where ${}_2F_1$ is the hypergeometric function.
  \end{ex}

  The unique minimizer to \eqref{eq:singularminimization} is therefore
  given by 
  $X(t)=\alpha(t)=X_{\mathrm{ref}}(t)$ for all
  $t\in[0,4]$. 
  Since  Example \eqref{ex:singular} has running cost $h$ and flux $f$
  with explicit time dependence, we introduce an extra state dimension,
  $s(t)=t$, as in Remark \ref{rem:timedep}.
  The Hamiltonian is then given by 
  \[
  H(x,s;\lambda_1,\lambda_2)=\frac{\lambda_1
  x}{\big((s-t_0)^2+\varepsilon^2\big)^{\beta/2}} - \frac{\lambda_1^2}{4\big((s-t_0)^2+\varepsilon^2\big)^{\beta}}+\lambda_2,
  \]
  where $\lambda_2$ is the dual corresponding to $s$.

\begin{figure}
\centering
  \includegraphics[width=0.7\textwidth]
  {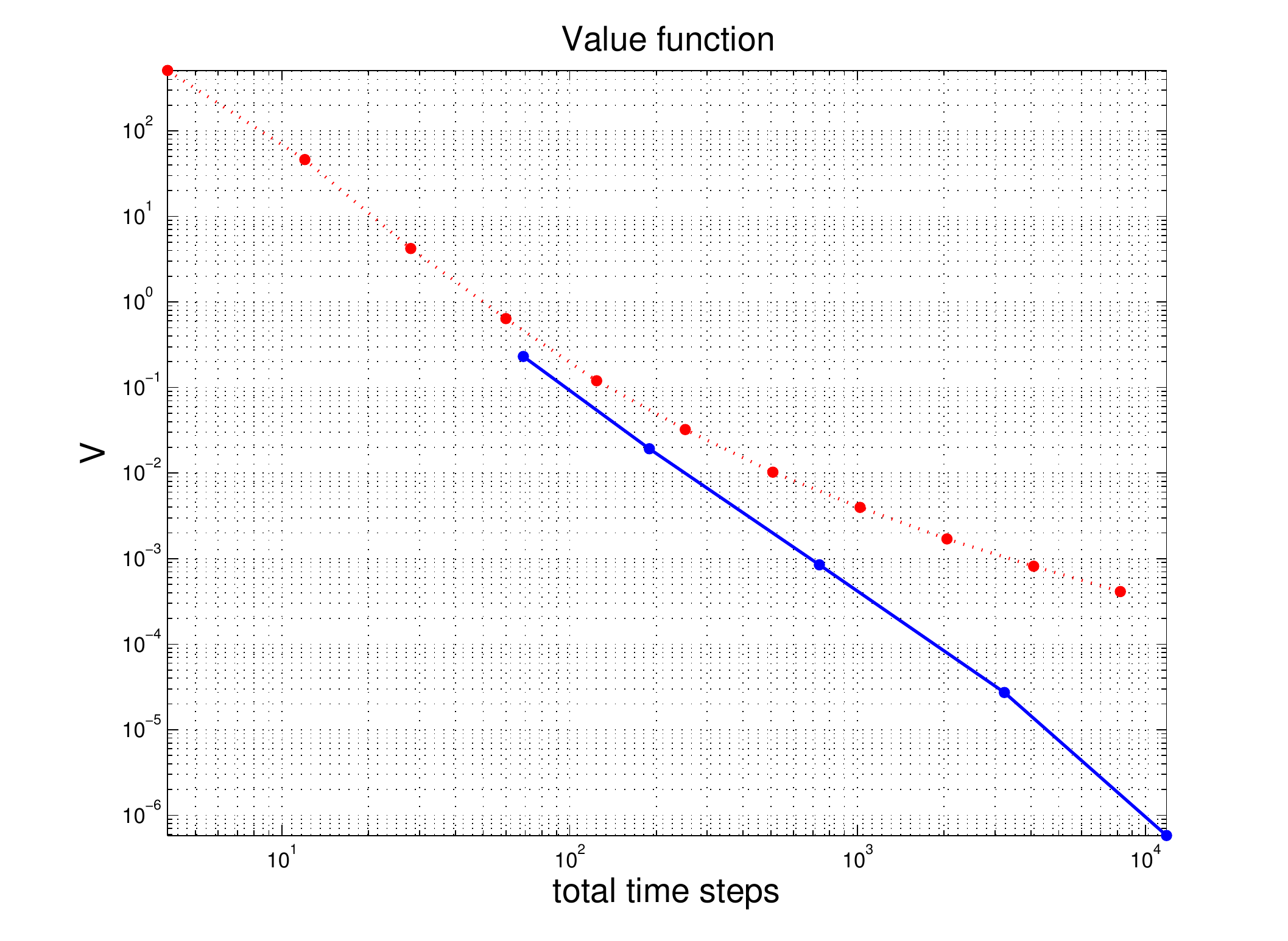}
  \caption{The minimum value of the functional in
    \eqref{eq:singularminimization} for the singular optimal control problem in Example \ref{ex:singular}, versus the cumulative number of time
    steps on all refinement levels for the adaptive algorithm (solid)
  and uniform time steps (dotted). Since the true value of
  \eqref{eq:singularminimization} is zero the graphs also indicate the
  respective errors. The regularization paramaters are
  $\varepsilon=10^{-10}$ and $\beta=3/4$.}
  \label{fig:SINGULAR_error}
\end{figure}

\begin{figure}
  \centering
  \includegraphics[width=0.7\textwidth]
  {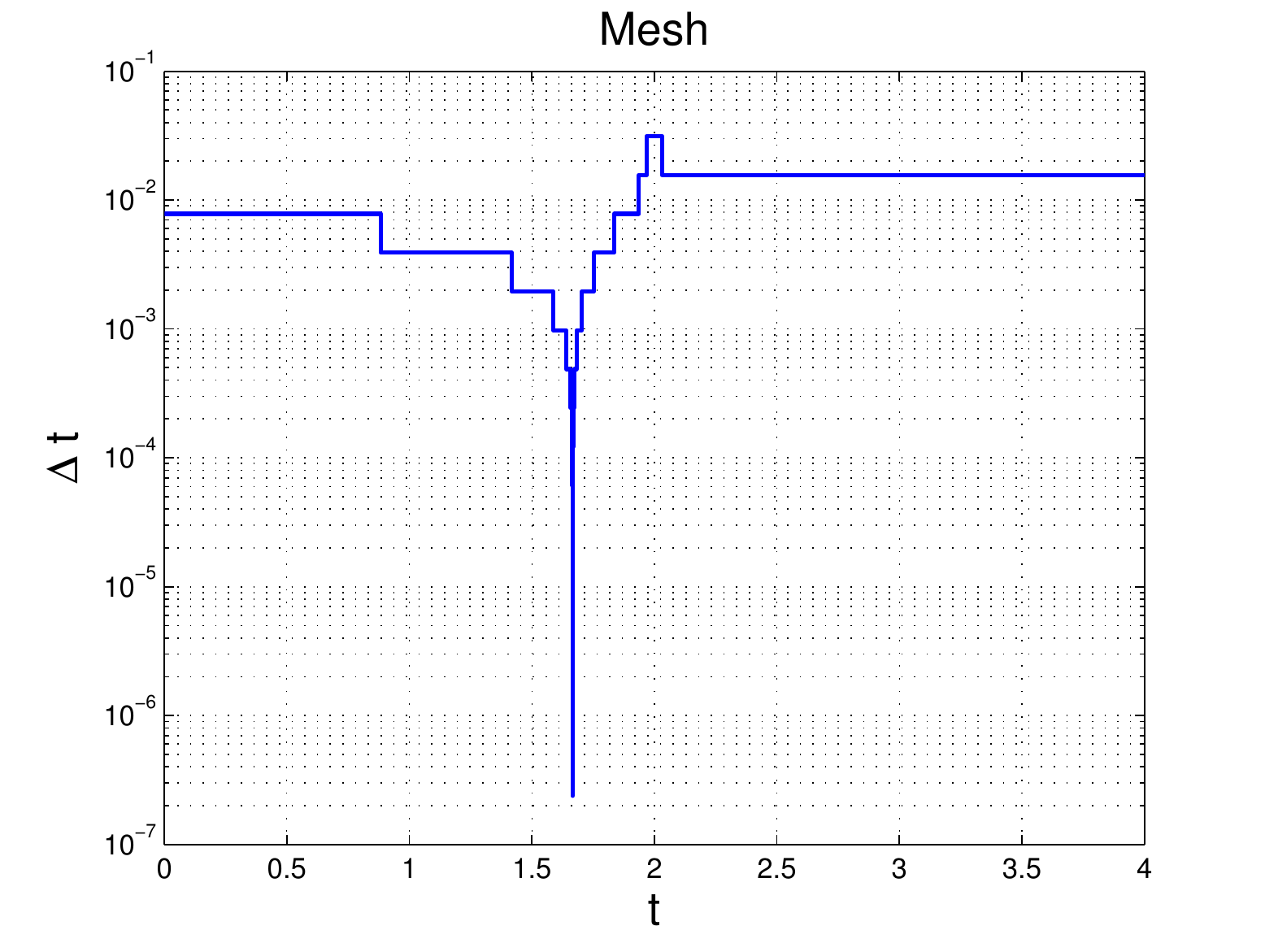}
  \caption{Mesh size versus time for the singular optimal control problem in Example \ref{ex:singular}. The regularization paramaters are 
  $\varepsilon=10^{-10}$ and $\beta=3/4$.}
  \label{fig:SINGULAR_mesh}
\end{figure}

  Although the Hamiltonian is a smooth function, the problem is a
  regularization of a  controlled ODE with a singularity, 
  \[
  X'(t)=\frac{\alpha(t)}{\abs{t-t_0}^\beta},
  \]
  and if the regularization parameter, $\varepsilon$, is small, the
  remainder term in Theorem \ref{thm:errorrepresentation} will be
  large unless the time steps are very small. 
  As the minimum value of the functional in
  \eqref{eq:singularminimization} is zero (attained for
  $\alpha=X=X_\mathrm{ref}$), it is immediately clear what the error in
  this functional is for a numerical simulation.   
  Figure
  \ref{fig:SINGULAR_error} shows 
  errors for adaptive and uniform time stepping versus the total number of
  time steps,
and Figure \ref{fig:SINGULAR_mesh} shows the dependence of the mesh size on
the time parameter.

\section{Conclusions}\label{sec:conclusions}
We have presented an a posteriori error representation 
for
optimal control problems with a 
bound for the remainder term.
With the error representation, it is possible to construct adaptive
algorithms, and we have presented and tested one such algorithm here. 
The error representation theorem assumes that the Hamiltonian
associated with the optimal control problem belongs to $C^2$. 
As many optimal control problems have Hamiltonians that are only
Lipschitz continuous, this is a serious restriction.
We have illustrated with a simple test example that 
$C^2$ smoothness may not be necessary. To justify this rigorously
remains an open problem.


\Appendix\section{Proof of Theorem \ref{thm:OptimEquivSP}}
\emph{Step 1.} We show here that there exist a constant $K$, and a continuous
function $S:[0,\infty)\rightarrow\Re$, such that
  $\lim_{s\rightarrow\infty} S(s)=\infty$, and 
\begin{equation}\label{eq:LSineq}
L(x,\beta)\geq (|\beta|-B|x|)_+ S\big((|\beta|-B|x|)_+\big)- K(1+|x|),
\end{equation}
where $y_+=\max\{y,0\}$.
We will show \eqref{eq:LSineq} with $K=\max\{\mu(0),A\}$ and $S$
defined by
\begin{equation*}
S(\xi)=\int_0\sp\xi\big|\big\{\chi:\mu'(\chi)\leq t, \chi\geq 0\big\}\big|\, \ud t/\xi.
\end{equation*}

We start by noting that the absolutely continuous (since it is convex) function $\mu$  can be
modified so that $\mu'>1$ almost everywhere  while \eqref{eq:Hmu}
still holds. We will henceforth assume that $\mu$ satisfies this
condition. 

By the bound on the Hamiltonian, $H$, and the definition of the
running cost, $L$, in \eqref{eq:LegendreL}, we have
\begin{equation*}
L(x,\beta)\geq \sup_{\lambda\in\Re\sp d}
\big\{-\beta\cdot\lambda-\mu(|\lambda|) - |x|(A+B|\lambda|)\big\}.
\end{equation*} 
By choosing $\lambda=-\chi\beta/|\beta|$, for $\chi\geq 0$, we have
\begin{equation*}
L(x,\beta) \geq \chi|\beta|-\mu(\chi)-|x|(A+B\chi)=:G_{x,\beta}(\chi).
\end{equation*}
Since $G_{x,\beta}(\cdot)$
is concave 
on $[0,\infty)$, at least one of the following
  alternatives must hold:
\renewcommand\theenumi{\Roman{enumi}}
\begin{enumerate}
\item \label{item1} $L(x,\beta)=\infty$.
\item \label{item2} $G_{x,\beta}$ is maximized at $\chi=0$.
\item \label{item3} $G_{x,\beta}$ is maximized at some $\chi\sp * \in (0,\infty )$.
\item \label{item4} $\sup_{0\leq \chi < \infty} G_{x,\beta}(\chi) =\lim_{\chi\rightarrow\infty}G_{x,\beta}(\chi)$.
\end{enumerate}

If alternative \ref{item1} holds, \eqref{eq:LSineq} is clearly satisfied with any $S$
and $K$. 
If alternative \ref{item2}  holds, then 
\begin{equation*}
L(x,\beta)\geq -\mu(0)-A\abs{x}.
\end{equation*}
Since $\chi=0$ maximizes $G_{x,\beta}$ and $\mu$ is convex it follows
that $S\big((\abs{\beta}-B\abs{x})_+\big)=0$. 
Hence \eqref{eq:LSineq} holds. 

If alternative \ref{item3} holds, we have 
\begin{equation*}
L(x,\beta)\geq (\abs\beta-B\abs x)\chi\sp *-\mu(\chi\sp *) -A\abs x.
\end{equation*}
Since $\mu$ is convex, it is absolutely continuous, and we have
\begin{equation*}
\mu(\chi\sp*)=\mu(0)+\int_0\sp{\chi\sp*}\mu'(\chi)\, \ud \chi.
\end{equation*}
Using  a layer cake representation (see \cite{Lieb-Loss}) of this
integral we get,
\begin{equation*}
\begin{split}
\int_0\sp{\chi\sp*}\mu'(\chi)\, \ud \chi &= \int_0\sp\infty
\big|\big\{\chi:\mu'(\chi)>t,\chi\in[0,\chi\sp*]\big\} \big| \, \ud t \\
&=\int_0\sp{\abs\beta-B\abs x}
\big|\big\{\chi:\mu'(\chi)>t,\chi\in[0,\chi\sp*]\big\} \big| \, \ud t,
\end{split}
\end{equation*}
where the absolute sign in the integrals denotes the Lebesgue measure,
and the last equality follows by the fact that $\mu'(\chi) \leq
\abs\beta-B \abs x$ for $\chi\in[0,\chi\sp*]$ since $\chi\sp*$ maximizes
$G_{x,\beta}(\chi)$. 
Since
\begin{equation*}
(\abs\beta-B\abs x)\chi\sp* = \int_0\sp{\abs\beta-B\abs x} \big|[0,\chi\sp*]\big|\, \ud t,
\end{equation*}
we have
\begin{equation*}
\begin{split}
(\abs\beta-B\abs x)\chi\sp*-\mu(\chi\sp*)&=-\mu(0)+\int_0\sp{\abs\beta-B\abs
    x}\big|\big\{ \chi: \mu'(\chi)\leq t,
    \chi\in[0,\chi\sp*]\big\}\big|\, \ud t\\
&=-\mu(0)+\int_0\sp{\abs\beta-B\abs x}\big|\big\{ \chi: \mu'(\chi)\leq t,
    \chi\geq 0\big\}\big|\, \ud t,
\end{split}
\end{equation*}
where the last inequality follows from the fact that $\mu'(\chi)\geq
\abs\beta -B\abs x$, when $\chi\geq \chi\sp*$. 
Since $\mu'$ is finite-valued almost everywhere we have
\begin{equation*}
\lim_{t\rightarrow\infty}\big|\big\{\chi:\mu'(\chi)\leq t, \chi\geq 0\big\}\big| =\infty,
\end{equation*}
and therefore $\lim_{s\rightarrow\infty}S(s)=\infty$.
Since $\mu'\geq 1$, the function $S$ is continuous.  
With $K=\max\{\mu(0),A\}$, \eqref{eq:LSineq} is satisfied.

If alternative \ref{item4} holds we can use that 
\begin{equation*}
L(x,\beta) \geq (\abs{\beta}-B\abs{x}-\varepsilon)\chi-\mu(\chi)-A\abs{x}=:G_{x,\beta}\sp{\varepsilon}(\chi)
\end{equation*}
for all $0\leq \chi<\infty$ and $\varepsilon>0$. For every
$\varepsilon>0$ the function $G_{x,\beta}$ is maximized at a
$\chi\sp{*}_\varepsilon\in[0,\infty)$. This gives, as the analysis for
  alternatives \ref{item2} and \ref{item3} shows, that 
\begin{equation*}
L(x,\beta)\geq (|\beta|-B|x|-\varepsilon)_+ S\big((|\beta|-B|x|-\varepsilon)_+\big)- K(1+|x|).
\end{equation*}
Since $\varepsilon$ could be chosen arbitrarily small and  positive
\eqref{eq:LSineq} follows.

\emph{Step 2.} We now show that for each time step $t_n$, there
exists a constant $K$, such that 
\begin{equation}\label{eq:ulowerbound}
\bar u(x,t_n)\geq -K(1+\abs x).
\end{equation} 
(The constant $K$ is allowed to depend on the time step $n$ and the
step length $\Delta t_n$.)
Assume \eqref{eq:ulowerbound} is satisfied at the time step
$t_{n+1}$. We will show that this implies that it is satisfied at
$t_n$ as well. 

The lower bound on $\bar u(\cdot,t_{n+1})$ and the lower bound on $L$
in \eqref{eq:LSineq}, together with dynamic programming gives
\begin{multline*}
\bar u(x,t_n)=\inf_{\beta\in\Re\sp d} \big(\Delta t_n L(x,\beta) +
\bar u(x+\Delta t_n \beta,t_{n+1})\big)\\
\geq 
\inf_{\beta\in\Re\sp
  d}\big(\Delta t_n(\abs\beta-B\abs x)_+ S\big((\abs\beta-B\abs x)_+\big) -
\tilde K - \tilde K\abs{x}-\tilde K
\abs{\beta}\big)=:\inf_{\beta\in\Re\sp d} J(x,\beta),
\end{multline*}
with a $\tilde K$ depending on $\Delta t_n$.
Since the function $S$ grows to infinity, there exists a $C\geq 0$, such
that $\xi\geq C$ implies $S(\xi)\geq \tilde K/\Delta t_n$. 
For such $\beta$  that satisfy $\abs\beta-B\abs x \geq C$ it therefore
holds that
\begin{equation*}
J(x,\beta)\geq \tilde K(\abs\beta-B\abs x) -\tilde K -\tilde K\abs x -\tilde K\abs\beta 
=-\tilde K - (\tilde K+\tilde K B)\abs x.
\end{equation*}
Since $S$ is continuous the function 
\begin{equation*}
\xi\mapsto \xi_+S(\xi_+)
\end{equation*}
attains a smallest value $D$ on the set 
$\{\xi\in\Re^d:\abs\xi\leq C\}$.
For every $\beta$ satisfying $\abs\beta-B\abs x\leq C$ we therefore have
\begin{equation*}
J(x,\beta)\geq D\Delta t_n-\tilde K-\tilde K\abs x-\tilde K\abs\beta
\geq D\Delta t_n-\tilde K-\tilde K C - (\tilde K+\tilde K B)\abs x.
\end{equation*}
With $\bar K=\max\{\tilde K + \tilde K B, \tilde K + \tilde K C -
D\Delta t_n\}$, and hence independent of $x$, we have
\begin{equation*}
\bar u(x,t_n) \geq - \bar K (1+\abs x).
\end{equation*}
Since $\bar u(\cdot,t_N)$ satisfies \eqref{eq:ulowerbound} with $K=k$,
by the lower bound on $g$, induction backwards in time shows that
\eqref{eq:ulowerbound} holds for all $n\leq N$, with different
constants $K$.


\emph{Step 3.} Assume that $\bar u(\cdot,t_{n+1})$ is locally
semiconcave. It is then also continuous (even locally Lipschitz
continuous, see e.g.\ \cite{Cannarsa-Sinestrari}). 
Since the Hamiltonian, $H$, is finite-valued everywhere, $L(x,\cdot)$
is lower semicontinuous, for every $x\in\Re\sp d$, see \cite{Clarke}.
Let
$\{\beta_i\}_{i=1}\sp\infty$ be a sequence of controls such that
\begin{equation*}
\lim_{i\rightarrow\infty} \Delta t_n L(X_n,\beta_i) + \bar u(X_n+\Delta
t\beta_i,t_{n+1}) \rightarrow \bar u(X_n,t_n).
\end{equation*}
By the lower bounds \eqref{eq:LSineq} and \eqref{eq:ulowerbound} for the functions
$L$ and $\bar u(\cdot,t_{n+1})$, proved in steps 1 and 2,  it follows that
the sequence $\{\beta_i\}_{i=1}\sp\infty$ is contained in a compact
set in $\Re\sp d$. It therefore contains a convergent subsequence 
\begin{equation*}
\beta_{i_j}\rightarrow\beta_n.
\end{equation*}
Since $\bar u(\cdot,t_{n+1})$ is continuous, and $L(X_n,\cdot)$ is
lower semicontinuous, we have that 
\begin{equation*}
\bar u(X_n,t_n)=\Delta t_n L(X_n,\beta_n)+\bar u(X_n+\Delta t_n\beta_n,t_{n+1}).
\end{equation*}

\emph{Step 4.} Assume that $\bar u(\cdot,t_{n+1})$ is locally semiconcave, and
that $\lambda_{n+1}$ is an element in $D\sp +\bar u(X_{n+1},t_{n+1})$, where
$X_{n+1}=X_n+\Delta t_n\beta_n$, and $\beta_n$ is the minimizer from
the previous step. We will
show that this implies
that 
\begin{equation}\label{eq:HLconnection}
\lambda_{n+1}\cdot\beta_n + L(X_n,\beta_n)=H(X_n, \lambda_{n+1}).\end{equation}

Consider a closed unit ball $B$ centered at $\beta_n$.
By the local semiconcavity of $\bar u(\cdot,t_{n+1})$ we have that there exists an $\omega:\Re\sp+
\rightarrow \Re\sp+$, such that $\lim_{\rho\rightarrow
  0\sp+}\omega(\rho)=0$, and
\begin{equation}\label{eq:onestepvalue}
\bar u(X_n+\Delta t_n \beta,t_{n+1}) \leq \bar u(X_{n+1},t_{n+1}) +
\Delta t_n \lambda_{n+1}\cdot(\beta-\beta_n) + \abs{\beta-\beta_n}\omega(\abs{\beta-\beta_n}),
\end{equation}
for all $\beta$ in $B$, see \cite{Cannarsa-Sinestrari}.
Since we know that the function 
\begin{equation*}
\beta\mapsto \bar u(X_n+\Delta t_n \beta,t_{n+1})+\Delta t_n L(X_n,\beta)
\end{equation*}
is minimized for $\beta=\beta_n$, the semiconcavity of $\bar u$ in 
\eqref{eq:onestepvalue} implies that  the function 
\begin{equation}\label{eq:funcalmostvalue}
\beta\mapsto \Delta t_n\lambda_{n+1}\cdot\beta
+\abs{\beta-\beta_n}\omega(\abs{\beta-\beta_n})  +\Delta t_n L(X_n,\beta)
\end{equation}
is also minimized on $B$ for $\beta=\beta_n$ (and therefore by the
convexity of $L(X_n,\cdot)$ also minimized on $\Re^d$). We will prove that 
the function 
\begin{equation}\label{eq:funcaffine}
\beta\mapsto \lambda_{n+1}\cdot\beta + L(X_n,\beta)
\end{equation}
is minimized for $\beta=\beta_n$. Let us assume that this is false,
so that there exists an $\beta\sp *\in\Re\sp d$, and an $\varepsilon > 0$, such that
\begin{equation}\label{eq:contradictionassumption}
\lambda_{n+1}\cdot\beta_n+L(X_n,\beta_n) -
\lambda_{n+1}\cdot\beta\sp*-L(X_n,\beta\sp *) \geq \varepsilon.
\end{equation}
Let $\xi\in[0,1]$, and $\hat\beta=\xi\beta\sp *
+(1-\xi)\beta_n$. Insert $\hat \beta$ into the function in
\eqref{eq:funcalmostvalue}:
\begin{equation*}
\begin{split}
&\Delta
t\lambda_{n+1}\cdot\hat\beta+|\hat\beta-\beta_n|\omega(|\hat\beta-\beta_n|)
+\Delta t_n L(X_n,\hat\beta)\\
&=
\Delta
t(\xi\lambda_{n+1}\cdot\beta\sp*+(1-\xi)\lambda_{n+1}\cdot\beta_n)+\xi\abs{\beta\sp*-\beta_n}\omega(\xi\abs{\beta\sp*-\beta_n})\\
&\qquad +\Delta t_n L(X_n,\xi\beta\sp *+(1-\xi)\beta_n)\\ 
&\leq \Delta
t(\xi\lambda_{n+1}\cdot\beta\sp*+(1-\xi)\lambda_{n+1}\cdot\beta_n)
+\xi\abs{\beta\sp*-\beta_n}\omega(\xi\abs{\beta\sp*-\beta_n})\\ 
&\qquad +\Delta t_n\xi L(X_n,\beta\sp *)
+\Delta t_n (1-\xi)L(X_n,\beta_n)\\
&\leq \Delta t_n(\lambda_{n+1}\cdot\beta_n + L(X_n,\beta_n))
+\xi\abs{\beta\sp*-\beta_n}\omega(\xi\abs{\beta\sp*-\beta_n})
-\Delta t_n \xi\varepsilon \\
&\qquad<\Delta t_n(\lambda_{n+1}\cdot\beta_n + L(X_n,\beta_n)),
\end{split}
\end{equation*}
for some small positive  number $\xi$. This contradicts the fact that
$\beta_n$ is a minimizer to the function in
\eqref{eq:funcalmostvalue}. Hence we have shown that the function in
\eqref{eq:funcaffine} is minimized at $\beta_n$. By the relation
between $L$ and $H$ in \eqref{eq:LegendreH} our claim
\eqref{eq:HLconnection} follows.

\emph{Step 5.} From the result in step 4, equation
\eqref{eq:HLconnection}, and the definition of the running cost $L$ in
\eqref{eq:LegendreL} it follows that
$\beta_n=H_\lambda(X_n,\lambda_{n+1})$, for if this equation did not
hold, then $\lambda_{n+1}$ could not be the maximizer of 
$-\beta_n\cdot\lambda+H(X_n,\lambda)$.


\emph{Step 6} We now show that under the assumption that $\bar
u(\cdot,t_{n+1})$ is locally semiconcave, then for each $F>0$ there
exists a $G>0$, such that 
\begin{equation}\label{eq:alphaboundimpl}
\abs{x}\leq F \implies \abs{\beta_x}\leq G,
\end{equation}
where $\beta_x$ is any optimal control at position $(x,t_n)$, i.e.\ 
$\bar u(x,t_n)=\bar u(x+\Delta t_n \beta_x,t_{n+1})+\Delta t_n L(x,\beta_x)$.
Step 5 proved that an optimal control is given by
$\beta_n=H_\lambda(X_n,\lambda_{n+1})$, so that 
\begin{equation*}
\bar u (0,t_n)=\bar u\big(\Delta t_n H_\lambda(0,p),t_{n+1}\big)+\Delta
t L\big(0,H_\lambda(0,p)\big),
\end{equation*}
where $p$ is an element in $D\sp + \bar u(\Delta t_n
\beta_0,t_{n+1})$. Let us now consider the control
$H_\lambda(x,p)$. Since this control is not necessarily optimal except
at $(0,t_n)$, we have
\begin{equation*}
\bar u (x,t_n)\leq \bar u \big(x+\Delta t_n
H_\lambda(x,p),t_{n+1}\big)+\Delta t_n L\big(x,H_\lambda(x,p)\big).
\end{equation*}
Since $\bar u(\cdot,t_{n+1})$ is locally semiconcave it is also
locally Lipschitz continuous (see \cite{Cannarsa-Sinestrari}).  
By the definition of $L$ in \eqref{eq:LegendreL} it follows that 
\begin{equation*}
L\big(x,H_\lambda(x,p)\big)=-H_\lambda(x,p)\cdot p+H(x,p).
\end{equation*} 
Since both $H(\cdot,p)$ and $H_\lambda(\cdot,p)$ are locally Lipschitz continuous
by assumption it follows that there exists a constant $E>0$ such that
\begin{equation}\label{eq:uliplower}
\bar u(x,t_n)-\bar u (0,t_n)\leq E,
\end{equation}
for all $\abs{x}\leq F$.

The inequalities \eqref{eq:LSineq} from step 1 and
\eqref{eq:ulowerbound} from step 2, together with \eqref{eq:uliplower}
give \eqref{eq:alphaboundimpl}.

\emph{Step 7.} In this step, we  show that if $\bar
u(\cdot,t_{n+1})$ is locally semiconcave, 
then so is $\bar u(\cdot,t_n)$. 
Furthermore, if $\beta_x$ is an optimal control at $(x,t_n)$, and $p$
is an element in $D\sp + \bar u(x+\Delta t_n \beta_x,t_n)$, then 
\begin{equation*}
p+\Delta t_n H_x(x,p) \in D\sp +\bar u(x,t_n).
\end{equation*}
We denote by $B_r$ the closed ball centered at the origin with
radius $r$.
In order to prove that $\bar u(\cdot,t_n)$ is locally semiconcave it is
enough to show that it is semiconcave on 
$B_r$, where $r$ is any positive radius. To accomplish this we will use the result
from step 6. We therefore take the radius $r=F$, which according to
step 6 can be taken arbitrarily large.

In step 3 we  showed that an optimal control $\beta_x$ exists at
every point  $x\in\Re\sp d$ at time $t_n$, under the assumption that
$\bar u(\cdot,t_{n+1})$ is locally semiconcave. 
In step 6 we showed that given any radius $F$ and $\abs{x}\leq F$,
there exists a constant $G$ such that all optimal controls must satisfy
$\abs{\beta_x}\leq G$.

A locally semiconcave function from $\Re\sp d$ to $\Re$ is locally
Lipschitz continuous (see \cite{Cannarsa-Sinestrari}). 
Hence, for every $x\in B_{F+G\Delta t_n}$, and every $p\in D\sp+\bar u(x,t_{n+1})$,
we have $\abs p \leq E$, for some constant $E$. By continuity, there
exists some constant $J$ such that 
$\abs{H_\lambda}\leq J$ on $B_F\times B_E$.

Let $R:=\max\{F+G\Delta t_n,F+J\Delta t_n\}$.
By the assumed local semiconcavity of $\bar u(\cdot,t_{n+1})$ we have
that there exists an $\omega:\Re\sp+\rightarrow\Re\sp+$, such that 
$\lim_{\rho\rightarrow 0}\omega(\rho)=0$, and
\begin{equation*}
\bar u (x,t_{n+1}) \leq \bar u(z,t_{n+1}) +p\cdot(x-z)+\abs{x-z}\omega(\abs{x-z}),
\end{equation*}
for all $x$ and $z$ in 
$B_R$, and $p$ in $D\sp +\bar u(z,t_{n+1})$, see
\cite{Cannarsa-Sinestrari}. We take $\omega$ to be non-decreasing,
which is clearly possible.
Let us now consider the control
$H_\lambda(x,p)$, 
where $p\in D\sp + \bar u(y+\Delta t_n\beta_y,t_{n+1})$, and $\beta_y$
is an optimal control at the point $y\in B_F$
($\beta_y=H_\lambda(y,p)$ according to step 5).
Since this
control is not necessarily optimal except at $(y,t_n)$, we have
\begin{multline}\label{eq:semicon1}
\bar u (x,t_n) \leq \bar u\big(x+\Delta t_n
H_\lambda(x,p),t_{n+1}\big) +\Delta t_n
L\big(x,H_\lambda(x,p)\big) \\
\leq \bar u(y+\Delta t_n\beta_y,t_{n+1})+p\cdot\big(x+\Delta t_n
H_\lambda(x,p)-(y+\Delta t_n\beta_y)\big)+ \Delta t_n
L\big(x,H_\lambda(x,p)\big)\\
+ |x+\Delta t_n
H_\lambda(x,p)-(y+\Delta t_n\beta_y)|\omega(|x+\Delta t_n
H_\lambda(x,p)-(y+\Delta t_n\beta_y)|). 
\end{multline}
By the bound on $\abs{H_\lambda}$,
this inequality holds for every $x$ and $y$ in $B_F$. 
By the definition of $L$ in \eqref{eq:LegendreL} it follows that 
\begin{equation}\label{eq:LHrel}
L\big(x,H_\lambda(x,p)\big)=-H_\lambda(x,p)\cdot p+H(x,p).
\end{equation} 
With this fact in \eqref{eq:semicon1}, and using that $\beta_y=H_\lambda(y,p)$, we have
\begin{multline}\label{eq:semicon2}
\bar u (x,t_n)
\leq \bar u(y+\Delta t_nH_\lambda(y,p),t_{n+1})+p\cdot(x-(y+\Delta t_nH_\lambda(y,p)))+ \Delta t_n H(x,p)\\
+ |x+\Delta t_n
H_\lambda(x,p)-(y+\Delta t_nH_\lambda(y,p))|\omega(|x+\Delta t_n
H_\lambda(x,p)-(y+\Delta t_nH_\lambda(y,p))|). 
\end{multline}
By the fact that 
$H_\lambda(\cdot,p)$ is locally Lipschitz continuous,
\begin{equation}\label{eq:semicon3}
|x-y + \Delta
 t\big(H_\lambda(x,p)-H_\lambda(y,p) \big)| 
\leq K|x-y|,
\end{equation}
for all $x$ and $y$ in $B_F$, and some 
constant $K$. 
We
also need the fact that
\begin{equation}\label{eq:semicon4}
\bar u(y,t_n)=\bar u(y+\Delta t_nH_\lambda(y,p),t_{n+1})+\Delta t_n L(y,H_\lambda(y,p)).
\end{equation}
We insert the results\eqref{eq:LHrel}, \eqref{eq:semicon3}, and \eqref{eq:semicon4},  into
\eqref{eq:semicon2} to get 
\begin{equation}\label{eq:semiconA}
\begin{split}
&\bar u(x,t_n)\\ 
&\leq \bar u(y,t_n)+p\cdot(x-y)+\Delta t_n 
\big(H(x,p)-H(y,p)\big) +K\abs{x-y}\omega(K\abs{x-y}) \\
&\leq  \bar u(y,t_n)+\big(p+\Delta t_n
H_x(y,p)\big)\cdot(x-y) +\abs{x-y}\tilde\omega(\abs{x-y})
\end{split}
\end{equation}
where
\begin{equation*}
\tilde\omega(\rho)=K\omega(K\rho)+\max\{\abs{H_x(z,q)-H_x(y,q)}:
\abs{z-y}\leq \rho,\ (z,y) \in B_F\times B_F\},
\end{equation*}
and $\lim_{\rho\rightarrow 0\sp +}\tilde\omega(\rho)=0$, since
$H_x$ is assumed to be continuous. 

We will now use equation \eqref{eq:semiconA} to show that $\bar
u(\cdot,t_n)$ is semiconcave on $B_F$. Let $x$ and $z$ be any elements
in $B_F$, and let $y=wx+(1-w)z$, where $w\in[0,1]$. As before, $p$ is
an element in $D\sp + \bar u(y+\Delta t_n\beta_y,t_{n+1})$.
The inequality in \eqref{eq:semiconA} with this choice of $y$ gives
\begin{multline}\label{eq:semiconB}
\bar u(x,t_n) \leq \bar u(wx+(1-w)z,t_n)\\ 
+(1-w)\big(p+\Delta t_n
H_x(wx+(1-w)z,p)\big)\cdot(x-z)+(1-w)\abs{x-z}\tilde\omega((1-w)\abs{x-z}),
\end{multline}
and with $x$ exchanged by $z$,
\begin{multline}\label{eq:semiconC}
\bar u(z,t_n) \leq \bar u(wx+(1-w)z,t_n)\\ 
+w\big(p+\Delta t_n
H_x(wx+(1-w)z,p)\big)\cdot(z-x)+w\abs{x-z}\tilde\omega(w\abs{x-z}).
\end{multline}
We multiply  \eqref{eq:semiconB} by $w$, and 
\eqref{eq:semiconC} by $1-w$, and add the resulting equations, to
get
\begin{equation*}
\begin{split}
&w\bar u(x,t_n)+(1-w)\bar u(z,t_n)\\ 
&\leq \bar u(wx+(1-w)z,t_n) + w(1-w)
\abs{x-z}\big(\tilde\omega((1-w)\abs{x-z})+\tilde\omega(w\abs{x-z})
\big)\\
&\leq \bar u(wx+(1-w)z,t_n) + w(1-w)
\abs{x-z}\hat\omega(\abs{x-z}),
\end{split}
\end{equation*}
if we let
\begin{equation*}
\hat\omega(\rho):=2\tilde\omega(\rho).
\end{equation*}
Since $x$ and $z$ can be any points in $B_F$, this shows 
that $\bar u(\cdot,t_n)$ is locally semiconcave. 

By \eqref{eq:semiconA} it also  follows that 
\begin{equation*}
p+\Delta t_n H_x(y,p)\in D\sp +\bar u(y,t_n).
\end{equation*}

\emph{Step 8.} 
Since $\bar u(x,T)=g(x)$, which is locally semiconcave, 
step 7  and  
induction backwards in time shows that $\bar u(\cdot,t_n)$ is locally
semiconcave for all $n$. 
In step 3 we showed that optimal controls exist at every position in
$\Re\sp d$ at time $t_n$, provided $\bar u(\cdot,t_{n+1})$ is locally
semiconcave. 
Hence
there exists a minimizer 
$(\beta_m,\ldots,\beta_{N-1})$ to the discrete minimization
functional $J_{(y,t_m)}$ in \eqref{eq:J}, for every $y\in\Re\sp d$
and $0\leq m\leq N$. Let $(X_m,\ldots,X_N)$ be a corresponding solution
to \eqref{eq:xalpharel}, and $\lambda_N$ an element in $D\sp +
g(X_N)$. From steps 5 and 7, we have that
$\beta_{N-1}=H_\lambda(X_{N-1},\lambda_N)$, and
$\lambda_{N-1}:=\lambda_N+\Delta t_{N-1} H_x(X_{N-1},\lambda_N)\in D\sp
+\bar u(X_{N-1},t_{N-1})$. Induction backwards in time shows that
there exists a dual path $\lambda_n$, $n=m,\ldots,N-1$, such that it
together with $X_n$, $n=m,\ldots,N$, satisfies the discretized
Hamiltonian system \eqref{eq:SymplPontr}.


\bibliographystyle{siam}
\bibliography{references}

\end{document}